\documentclass[10pt]{article}
\usepackage[cp1252]{inputenc}
\usepackage[T1]{fontenc}
\usepackage{lmodern}
\usepackage[a4paper]{geometry}
\usepackage[english]{babel}
\usepackage{pifont}
\usepackage{extsizes}
\usepackage{xcolor}
\usepackage{fancybox}
\usepackage{graphicx}
\usepackage{sectsty}
\usepackage{enumerate}
\usepackage{multicol}
\usepackage[french]{varioref}
\usepackage[fleqn]{amsmath}
\usepackage{mathtools}
\usepackage{amssymb}
\usepackage{mathbbol}
\usepackage{mathrsfs}
\usepackage{bm}
\usepackage{amsthm}
\usepackage{blkarray}
\theoremstyle{break}
\newtheorem{Th}{Theorem}
\newtheorem{Prop}[Th]{Proposition}

\newtheorem{Le}{Lemma}
\newtheorem{Cor}{Corollary}
\theoremstyle{definition}
\newtheorem{Rk}{Remark}
\newtheorem{Ex}{Example}

\newcommand{\mm}{\mathfrak{m}}

\newcommand{\PP}{\mathbb{P}}
\newcommand{\EE}{\mathbb{E}}
\newcommand{\RR}{\mathbf{R}}
\newcommand{\NN}{\mathbb{N}}
\newcommand{\CC}{\mathfrak{C}_{\rho}}

\newcommand{\cell}{\mathcal{C}}
\newcommand{\conv}[2][n]{\underset{#1\rightarrow #2}{\longrightarrow}}
\newcommand{\EEE}[1]{\mathbb{E}\left[#1 \right]}
\newcommand{\PPP}[1]{\mathbb{P}\left(#1 \right)}
\newenvironment{prooft}[1]{\vskip 2mm\noindent {\bf Proof of #1.}}
                    {\hfill $\square$ \vskip 2mm \noindent}

\begin{document}
\pagenumbering{arabic}
\pagestyle{plain}
\author{ \textbf{Nicolas Chenavier}\footnote{Postal address: Université de Rouen, LMRS, avenue de l'Université, BP 12
76801 Saint-Etienne-du-Rouvray cedex, France. E-mail: nicolas.chenavier@etu.univ-rouen.fr}}
\title{\textbf{A general study of extremes of stationary tessellations with applications}}
\maketitle

\begin{abstract}
Let $\mathfrak{m}$ be a random tessellation in $\mathbf{R}^d$ observed in a bounded Borel subset $W$ and $f(\cdot)$ be a measurable function defined on the set of convex bodies. To each cell $C$ of $\mathfrak{m}$ we associate a point $z(C)$ which is the nucleus of $C$. Applying $f(\cdot)$ to all the cells of $\mathfrak{m}$, we investigate the order statistics of $f(C)$ over all cells $C\in\mathfrak{m}$  with nucleus  in $\mathbf{W}_{\rho}=\rho^{1/d}W$ when $\rho$ goes to infinity. Under a strong mixing property and a local condition on $\mathfrak{m}$ and $f(\cdot)$, we show a general theorem which reduces the study of the order statistics to the random variable $f(\mathscr{C})$ where $\mathscr{C}$ is the typical cell of $\mathfrak{m}$. The proof is deduced from a Poisson approximation on a dependency graph via the Chen-Stein method. We obtain that the point process $\left\{(\rho^{-1/d}z(C), a_{\rho}^{-1}(f(C)-b_{\rho})), C\in\mathfrak{m}, z(C)\in \mathbf{W}_{\rho}\right\}$, where $a_{\rho}>0$ and $b_{\rho}$ are two suitable functions depending on $\rho$, converges to a non-homogeneous Poisson point process. Several applications of the general theorem are derived in the particular setting of Poisson-Voronoi and Poisson-Delaunay tessellations and for different functions $f(\cdot)$ such as the inradius, the circumradius, the area, the volume of the Voronoi flower and the distance to the farthest neighbor. When the local condition does not hold and the normalized maximum converges, the asymptotic behaviour depends on two quantities that are the distribution function of $f(\mathscr{C})$ and a constant $\theta\in [0,1]$ which is the so-called extremal index.
\end{abstract}

  

\vspace{0.5cm}


\vspace{0.3cm}


\section{Introduction} 
\label{sectionintrochap2}
A tessellation of $\RR^d$, endowed with its natural norm $|\cdot|$, is a countable collection of compact subsets, called \textit{cells}, with disjoint interiors which subdivides the space and such that the number of cells intersecting any bounded subset of $\RR^d$ is finite. By a random  tessellation $\mm$, we mean a random variable defined on a hypothetical probability space $(\Omega, \mathcal{A}, \PP)$ with values in the set of tessellations of $\RR^d$ endowed with a specific $\sigma$-algebra induced by the Fell topology. It is said to be stationary if its distribution is invariant under translation of the cells. For a complete account on random tessellations, we refer to the books \cite{SW}, \cite{SKM}  and the survey \cite{Cal5}.

Given a fixed realization of $\mm$, we associate to each cell $C\in\mathfrak{m}$ in a deterministic way a point $z(C)$,  which is called the \textit{nucleus} of the cell, such that $z(C+x)=z(C)+x$ for all $x\in\RR^d$. To describe the mean behaviour of the tessellation, the notions of intensity and typical cell are introduced as follows. Let $B$ be a Borel subset of $\RR^d$ such that $\lambda_d(B)\in (0,\infty)$ where $\lambda_d$ is the $d$-dimensional Lebesgue measure. The \textit{intensity} $\gamma$ of the tessellation is defined as 
\[\gamma = \frac{1}{\lambda_d(B)}\cdot\EEE{\#\{C\in \mm, z(C)\in B\}}\] and we assume that $\gamma\in (0,\infty)$. Since $\mm$ is stationary, $\gamma$ is independent of $B$ and we suppose, without loss of generality, that $\gamma=1$. The \textit{typical cell} $\cell$ is a random polytope such that the distribution is given by  \begin{equation}\label{campbell}\EE[f(\cell)] = \frac{1}{\lambda_d(B)}\cdot\EEE{\sum_{\underset{z(C)\in B}{C\in\mathfrak{m},} }f(C-z(C))}\end{equation} where $f:\mathcal{K}_d\rightarrow\RR$ is any bounded measurable function on the set of convex bodies $\mathcal{K}_d$ (endowed with the Hausdorff topology). 

We are interested in the following problem: only a part of the tessellation is observed in the window $\mathbf{W}_{\rho}=\rho^{1/d}W$ where $W$ is a bounded Borel subset of $\RR^d$, i.e. included in a cube $\mathbf{C}^{(W)}$, and such that $\lambda_d(W)\neq 0$. Let $f:\mathcal{K}_d\rightarrow\RR$ be a translation invariant measurable function, i.e. $f(C+x)=f(C)$ for all $C\in\mathcal{K}_d$ and $x\in\RR^d$. We denote by $M^{(r)}_{f,\mathbf{W}_{\rho}}$ the $r$-th order statistic of $f$ over the cells $C\in\mm$ such that $z(C)\in \mathbf{W}_{\rho}$. When $r=1$, the 1-st order statistic is denoted by $M_{f,\mathbf{W}_{\rho}}$ i.e.
\[M_{f,\mathbf{W}_{\rho}}=M^{(1)}_{f,\mathbf{W}_{\rho}}=\max_{\underset{z(C)\in \mathbf{W}_{\rho}}{C\in\mathfrak{m},}}f(C).\] In this paper, we investigate the limit behaviour of $M^{(r)}_{f,\mathbf{W}_{\rho}}$ when $\rho$ goes to infinity. 

The study of extremes describes the regularity of the tessellation. For instance, in finite element method, the quality of the approximation depends on some consistency measurements over the partition, see e.g. \cite{JGZ}. Another potential application field is statistics of point processes. The key idea would be to identify a point process from the extremes of a tessellation induced by the point process. 

To the best of our knowledge, one of the first works on extreme values in stochastic geometry is due to Penrose. In chapters 6,7 and 8 in \cite{Pr}, he investigates the maximum and minimum degrees of random geometric graphs. More recently, Schulte and Thäle \cite{ST} establish a theorem to derive the smallest values of a functional $f_k(x_1,\ldots, x_k)$ of $k$ points on a homogeneous Poisson point process. Nevertheless, their approach cannot be applied to our problem for several reasons: first, they consider a Poisson point process. Moreover, studying extremes of the tessellation requires to use functionals which depend on the whole point process of nuclei and not only on a fixed number of points. In this paper, we consider any function $f(\cdot)$ and we restrict our investigation to a certain kind of random tessellation satisfying a strong mixing property. We give a general theorem, with the rates of convergence, which is followed by numerous examples in the particular setting of Poisson-Voronoi and Poisson-Delaunay tessellations. This improves in particular some extremes that are investigated in \cite{CC}. Before stating our main theorems, we need some preliminaries which contain notations and conditions on the random tessellation.

\paragraph*{Preliminaries.}
Let $\mathbf{C}^{(W)}$ be a cube in $\RR^d$ containing $W$. We partition $\mathbf{C}_{\rho}^{(W)}=\rho^{1/d}\mathbf{C}^{(W)}$ by a set $V_{\rho}$ of $N_{\rho}$ sub-cubes of equal size with $N_{\rho}\conv[\rho]{\infty}\infty$. These sub-cubes are denoted by indices $\mathbf{i} = (i_1,\ldots, i_d)\in V_{\rho}$. Let us define a distance between sub-cubes $\mathbf{i}$ and $\mathbf{j}$ as
\[d(\mathbf{i}, \mathbf{j}) = \max_{1\leq r\leq d}\{|i_r-j_r|\}.\] Moreover, if $A$, $B$ are two sets of sub-cubes, we let
$d(A,B) = \min_{\mathbf{i}\in A, \mathbf{j}\in B}d(\mathbf{i},\mathbf{j})$. 
For each $\mathbf{i}\in V_{\rho}$, we denote by \[ M_{f,\mathbf{i}} =  \max_{\underset{z(C)\in\mathbf{i}\cap \mathbf{W}_{\rho}}{C\in\mm,}}f(C). \] When $\{C\in\mathfrak{m}, z(C)\in \mathbf{i}\cap \mathbf{W}_{\rho}\}$ is empty, we take $ M_{f,\mathbf{i}}=-\infty$. 

Let us consider a threshold $v_{\rho}$ that is a function depending on $\rho$. Studying the order statistics amounts to investigate the number of exceedance cells $U_{\rho}(v_{\rho})$ defined as
\begin{equation}\label{defUgamma1}U_{\rho}(v_{\rho}) = \sum_{\underset{z(C)\in \mathbf{W}_{\rho}}{C\in\mm,}}\mathbb{1}_{f(C)>v_{\rho}}.\end{equation}  Thanks to \eqref{campbell}, the mean of this random variable is
\begin{equation}\label{meanexceedance} \EEE{U_{\rho}(v_{\rho}) } = \lambda_d(\mathbf{W}_{\rho})\cdot \PPP{f(\cell)>v_{\rho}}.\end{equation} 
We assume the following condition which is referred as the typical cell property (TCP):

\bigskip

\textsc{Condition (TCP)}: \textit{the mean number of exceedance cells converges to a limit denoted by $\tau\geq 0$ i.e.}
\[ \lambda_d(\mathbf{W}_{\rho})\cdot \PPP{f(\cell)>v_{\rho}}\conv[\rho]{\infty}\tau.\]
Moreover, we denote by $G_1(\rho)$ the rate of convergence i.e.
\begin{equation}\label{defG1general}G_1(\rho) = \left|\lambda_d(\mathbf{W}_{\rho})\cdot\PP(f(\cell)>v_{\rho}) - \tau\right|. \end{equation}

We assume also a (global) condition of $R$-dependence associated to $\mm$ and $f$ which is referred as \textsc{Condition 1}.  

\bigskip
\textsc{Condition 1}: \textit{there exists an integer $R$  and an event $A_{\rho}$ with $\PPP{A_{\rho}}\conv[\rho]{\infty}1$ such that, conditional on $A_{\rho}$ , the $\sigma$-algebras $\sigma\{M_{f,\mathbf{i}}, \mathbf{i}\in A\}$ and $\sigma\{M_{f,\mathbf{i}}, \mathbf{i}\in B\}$ are independent when $d(A,B)>R$.} 
\bigskip

Finally, in order to present our first theorem, we introduce a second function defined as
\begin{equation}\label{defG2general}G_2(\rho) =N_{\rho}\EE\left[\sum_{\underset{z(C_1), z(C_2)\in \mathfrak{C}_{\rho}}{(C_1,C_2)_{\neq}\in\mathfrak{m}^2,} }\mathbb{1}_{f(C_1)>v_{\rho}, f(C_2)>v_{\rho}} \right]\end{equation} where \begin{equation}\label{defcubeC}\mathfrak{C}_{\rho} = \left[0,(2R+1)\cdot\lambda_d(\mathbf{W}_{\rho})^{1/d}N_{\rho}^{-1/d}\right]^d\end{equation} and where $(C_1,C_2)_{\neq}\in\mathfrak{m}^2$ means that $(C_1,C_2)$ is a couple of distinct cells of $\mm$. 

\paragraph{Order statistics}

We are now prepared to present our first theorem.

\begin{Th}
\label{realtyp}
Let $\mathfrak{m}$ be a stationary random tessellation of intensity $1$ such that  \textsc{Condition (TCP)}  and \textsc{Condition 1} hold. Then  
\[\left|\PP(M^{(r)}_{f,\mathbf{W}_{\rho}}\leq v_{\rho}) - e^{-\tau}\sum_{k=0}^{r-1} \frac{\tau^k}{k!}  \right| = O\left(N_{\rho}^{-1} + \PP(A_{\rho}^c) + G_1(\rho) + G_2(\rho) \right).\] 
where $\phi(\rho) = O(\psi(\rho))$ means that $\phi(\rho)/\psi(\rho)$ is bounded.

\end{Th} 
To derive useful applications, we assume a second condition on the random tessellation.

\bigskip
\textsc{Condition 2}: \textit{the function $G_2(\rho)$ converges to 0 as $\rho$ goes to infinity.}
\bigskip

This (local) condition means that with high probability two neighbor cells are not simultaneously exceedances. With this assumption, we obtain the following result:

\begin{Cor}
Let $\mathfrak{m}$ be a stationary random tessellation of intensity $1$ such that  \textsc{Condition (TCP)}  and \textsc{Conditions 1} and \textsc{2} hold. Then
\[\PP(M^{(r)}_{f,\mathbf{W}_{\rho}}\leq v_{\rho})\conv[\rho]{\infty}e^{-\tau}\sum_{k=0}^{r-1} \frac{\tau^k}{k!}.\] 
\end{Cor}

The rate of convergence is provided in Theorem \ref{realtyp}. Besides, Theorem \ref{realtyp} could be extended to more general models such as Boolean models and marked point processes. When the random tessellation is ergodic with respect to the group of tessellations of $\RR^d$, the order statistics 
are asymptotically independent of the choice of nuclei $z(\cdot)$. This will be the case for the examples that we deal with. Indeed, they only depend on the asymptotic behaviour of $G_1(\rho)$ given by \eqref{defG1general} and the typical cell $\cell$ itself does not depend on the set of nuclei  thanks to Wiener ergodic's theorem. Moreover, we notice that the order statistics do not depend on the shape of the window $W$. Actually, a method similar to Proposition 3 of \cite{CC} shows that the contribution of boundary cells is negligible. 

As mentioned above, \textsc{Conditions 1 } and \textsc{2} concern global and local properties of the tessellation respectively. In fact, there exists an analogy between \textsc{Conditions 1} and \textsc{2} and Conditions $D(u_n)$ and $D'(u_n)$ of Leadbetter \cite{L1} respectively. The general theory of extreme values deals with sequences \cite{HT} or random fields \cite{Ch}, \cite{LR}, see also the reference books \cite{HF} and  \cite{R}. Unfortunately, we are unable to apply it in our setting. Indeed, the set of random variables that we consider is not a discrete random field in a classical meaning. More precisely, the process $\{M_{f,\mathbf{i}}\}_{\mathbf{i}\in V_{\rho}}$ is a triangular array indexed by $\NN^d$ and the process $\{f(C_x)\}_{x\in\RR^d}$ is not a Gaussian continuous random field, where $C_x$ is the cell of the tessellation containing $x$. 

\paragraph{Point process of exceedances}

In practice, the threshold is often of the form $v_{\rho}=v_{\rho}(t)=a_{\rho}t+b_{\rho}$, $t\in\RR$ with $a_{\rho}>0$. In that case, we can be more specific about the joint distributions of the order statistics. Before stating our second theorem, we need some preliminaries. We denote by $\tau(t)\in [0,+\infty]$, $t\in\RR$, the limit of $\lambda_d(\mathbf{W}_{\rho})\cdot\PP(f(\cell)>v_{\rho}(t))$ and by ${}_*x=\inf\{t\in\RR, \tau(t)<\infty\}$ and $x^*=\sup\{t\in\RR, \tau(t)>0\}$ the lower and upper endpoints of $\tau(\cdot)$. Since $a_{\rho}$ is positive, the function $\tau(\cdot)$ is not increasing so that $\tau(\cdot)$ is finite on $({}_*x,x^*]$.

 Under \textsc{Conditions 1} and \textsc{2}, we consider the random collection
\[\Phi_{\rho} = \left\{\left(\rho^{-1/d}z(C), a_{\rho}^{-1}(f(C) - b_{\rho} ) \right), C\in\mm \text{ and } z(C)\in \mathbf{W}_{\rho} \right\}\subset W\times \RR.\] Moreover, we consider a Poisson point process $\Phi\subset W\times ({}_*x,x^*]$,  with intensity measure $\nu$ given by \[\nu(B\times (s,t])=\EEE{\#\Phi\cap (B\times (s,t])} = \frac{\lambda_d(B)}{\lambda_d(W)}\cdot (\tau(s)-\tau(t))\] for all Borel subset $B\subset W$ and all segment $(s,t]\subset ({}_*x,x^*]$. 
We then obtain the following limit theorem.

\begin{Th}
\label{PPTh}
Let $\mathfrak{m}$ be a stationary random tessellation of intensity $1$ such that  \textsc{Condition (TCP)}  and \textsc{Conditions 1} and \textsc{2} hold for each $v_{\rho} = v_{\rho}(t)=a_{\rho}t+b_{\rho}$, $t\in \RR$. Then the family of point processes $\Phi_{\rho}$ converges in distribution to the Poisson point process $\Phi$ i.e. for any Borel subset $\mathscr{B}_1,\ldots, \mathscr{B}_k\subset W\times ({}_*x,x^*]$ with $\nu(\partial \mathscr{B}_i)=0$ for all $i=1,2,\ldots, k$
\[\left(\#\Phi_{\rho}\cap \mathscr{B}_1,\ldots, \#\Phi_{\rho}\cap \mathscr{B}_k \right) \overset{\mathcal{D}}{\longrightarrow} \left(\#\Phi\cap \mathscr{B}_1,\ldots, \#\Phi\cap \mathscr{B}_k \right)\] where $\partial\mathscr{B}$ denotes the boundary of $\mathscr{B}\subset W\times ({}_*x,x^*]$. 
\end{Th}
This result suggests that the largest order statistics can be seen as points of a (non homogeneous) Poisson point process. Theorem \ref{PPTh} gives their joint distributions so that Theorem \ref{realtyp} is a particular case of the latter when $k=1$ and $B=W\times (t,\infty)$. For a wider panorama on results of the point process of exceedances associated to the extremes of a sequence of non independent random variables, we refer to chapter 5 in \cite{LLR}. When $W=\mathbf{C}^{(W)}=[0,1]^d$ and when $\tau(\cdot)$ is not constant, the function $\tau(\cdot)$ belongs to either the Fréchet, the Gumbel or the Weibull family. This fact is a rewriting of the proof of Theorem 4.1 in \cite{LR}. 

\paragraph{Extremal index}
When \textsc{Condition 2} does not hold, the exceedance locations can be divided into clusters and the order statistics cannot be investigated when $r\geq 2$. Yet, the behaviour of $M_{f,\mathbf{W}_{\rho}}$ can be deduced up to a constant according to the following proposition. For sake of simplicity, we assume in Proposition \ref{realtyp2} that $W=[0,1]^d$.

\begin{Prop}
\label{realtyp2}
Let $\mathfrak{m}$ be a stationary random tessellation of intensity $1$ such that \textsc{Condition 1} holds and let $W=[0,1]^d$. Let us assume that for all $\tau\geq 0$, there exists a deterministic function $v_{\rho}(\tau)$ depending on $\rho$ such that $\rho\cdot\PP(f(\cell)>v_{\rho}(\tau))$ converges to $\tau$ as $\rho$ goes to infinity. Then there exist constants $\theta, \theta', 0\leq \theta\leq\theta'\leq 1$ such that, for all $\tau \geq 0$, 
\[ \limsup_{\rho\rightarrow\infty}\PP(M_{f,\mathbf{W}_{\rho}}\leq v_{\rho}(\tau)) = e^{-\theta\tau} \text{ and }  \liminf_{\rho\rightarrow\infty}\PP(M_{f,\mathbf{W}_{\rho}}\leq v_{\rho}(\tau)) = e^{-\theta'\tau} .\] In particular, if $\PPP{M_{f,\mathbf{W}_{\rho}}\leq v_{\rho}(\tau)}$ converges, then $\theta = \theta'$ and \[\PPP{M_{f,\mathbf{W}_{\rho}}\leq v_{\rho}(\tau)} \conv[\rho]{\infty}e^{-\theta\tau}.\]
\end{Prop}
Proposition \ref{realtyp2} is similar to the result due to Leadbetter for stationary sequences of real random variables (see
Theorem 2.2 of \cite{L2}). Its proof relies notably on the adaptation to our setting of several arguments included in \cite{L2}.  According to Leadbetter, we say that the random tessellation $\mathfrak{m}$ has \textit{extremal index} $\theta$ if, for each $\tau\geq 0$, $\rho\cdot\PPP{f(\cell)>v_{\rho}(\tau)}\conv[\rho]{\infty}\tau$ and $\PPP{M_{f,\mathbf{W}_{\rho}}\leq v_{\rho}(\tau)}\conv[\rho]{\infty}e^{-\theta\tau}$. In a future work, we hope to develop a general method to estimate the extremal index. 

\bigskip

The paper is organized as follows. In section \ref{sectionrealtyp}, we show how to reduce our problem to the study of extreme values on a dependency graph. We use a result of \cite{AGG} to derive an estimation of exceedances by a Poisson distribution. We then deduce Theorems \ref{realtyp} and \ref{PPTh} from a discretization of $W$ into sub-cubes. Sections \ref{sectionPDT}, \ref{sectionPVT} and \ref{sectionGP} are devoted to numerous applications on Delaunay and Voronoi random tessellations. We investigate the asymptotic behaviours with the rates of convergence of :
\begin{itemize}
\item the minimum of circumradii of a Poisson-Delaunay tessellation in any dimension and the maximum and minimum of the areas in the planar case (section \ref{sectionPDT}),
\item the minimum of distances to the farthest neighboring nucleus and the minimum of the volume of flowers for a Poisson-Voronoi tessellation (section \ref{sectionPVT}),
\item the maximum of inradii for a Voronoi tessellation induced by a Gauss-Poisson process (section \ref{sectionGP}).
\end{itemize}
For each tessellation and each characteristic, we need to find a suitable threshold $v_{\rho}$ and to check \textsc{Condition 2} which requires some delicate geometric estimates. In the last section, we prove Proposition \ref{realtyp2} and we give two examples where the extremal index differs from 1.

In the rest of the paper, $c$ or $c'$ denotes a generic constant which does not depend on $\rho$ but may depend on other quantities. The term $v_{\rho}=v_{\rho}(t)$ denotes a generic function of $t$, depending on $\rho$, which is specified in  sections \ref{sectionPDT}, \ref{sectionPVT} and \ref{sectionGP}.

\section{Proofs of Theorems \ref{realtyp} and \ref{PPTh} }
\label{sectionrealtyp}

\subsection{Extreme values on a dependency graph and proof of Theorem \ref{realtyp}}
We first outline the methodology of the proof of Theorem \ref{realtyp} with some additional notations. A classical method in extreme value theory is to investigate the exceedances. We consider two random variables that are  the number of exceedance cells $U_{\rho}(v_{\rho})$, introduced in \eqref{defUgamma1}, and the number of exceedance cubes $U'_{V_{\rho}}(v_{\rho})$ defined as
\begin{equation}\label{defUgammaUV}U_{\rho}(v_{\rho}) = \sum_{\underset{z(C)\in \mathbf{W}_{\rho}}{C\in\mm,}}\mathbb{1}_{f(C)>v_{\rho}} \text{ and } U'_{V_{\rho}}(v_{\rho}) = \sum_{\mathbf{i}\in V_{\rho}}\mathbb{1}_{M_{f,\mathbf{i}}>v_{\rho}} \end{equation} where $V_{\rho}$ and $M_{f,\mathbf{i}}$ are introduced in the preliminaries. We denote by $\mu_{\rho}$ the mean of $U'_{V_{\rho}}(v_{\rho})$ i.e .\begin{equation}\label{deftaugamma}\mu_{\rho} = \EEE{U'_{V_{\rho}}(v_{\rho})} = \sum_{\mathbf{i}\in V_{\rho}}\PPP{M_{f,\mathbf{i}}>v_{\rho}}.\end{equation} The proof of Theorem \ref{realtyp} can be displayed as the three following results.

\begin{Le}
\label{Thle1}
With the same assumptions as in Theorem \ref{realtyp}, we get for all $r\in\NN^*$
\begin{equation}\label{majThle1}\left|\PPP{U_{\rho}(v_{\rho})\leq r-1} - \PPP{U'_{V_{\rho}}(v_{\rho})\leq r-1} \right| \leq 2\cdot G_2(\rho)\end{equation}
\end{Le}
The above lemma is a consequence of \textsc{Condition 2}. 

\begin{Le}
\label{Thle2}
Let $\mu_{\rho}$ be as in \eqref{deftaugamma}. With the same assumptions as in Theorem \ref{realtyp}, we get for all $r\in\NN^*$
\begin{equation}\label{majUVgamma}\left|\PP(U'_{V_{\rho}}(v_{\rho})\leq r-1) - e^{-\mu_{\rho}}\sum_{k=0}^{r-1}\frac{\mu_{\rho}^k}{k!} \right| = O\left(N_{\rho}^{-1} + \PP(A_{\rho}^c) + G_2(\rho)  \right)\end{equation}
\end{Le}
The derivation of the latter constitutes the major part of the proof of Theorem \ref{realtyp}. It means that the number of exceedance cubes is approximately a Poisson random variable. The fundamental concept to prove this lemma is that of a  dependency graph. We first establish a Poisson approximation on the number of exceedances on such graph and we show how we can reduce our problem to this graph. Finally, the following result gives an estimate of $\mu_{\rho}$. 

\begin{Le}
\label{Thle3}
Let $\mu_{\rho}$ as in \eqref{deftaugamma}. With the same assumptions as in Theorem \ref{realtyp}, we get
 \begin{equation}\label{meanestimateexp}|\mu_{\rho}-\tau| \leq G_1(\rho)+ G_2(\rho). \end{equation}
\end{Le}

\begin{prooft}{Theorem \ref{realtyp}}
Since $M_{f,\mathbf{W}_{\rho}}^{(r)}$ is lower than $v_{\rho}$ if and only if $U_{\rho}(v_{\rho})\leq r-1$, we deduce Theorem \ref{realtyp} from the three lemmas above and the fact that the function $x\mapsto e^{-x}\sum_{k=0}^{r-1}\frac{x^k}{k!}$ is Lipschitz. 
\end{prooft}

In the rest of the subsection, we proceed as follows:
\begin{enumerate}
\item We establish the Poisson approximation on a dependency graph (Proposition \ref{steingraph}) and we deduce from it  Lemma \ref{Thle2}. The key idea is to apply \textsc{Condition (TCP)} and \textsc{Condition 1}.
\item We prove Lemmas \ref{Thle1} and \ref{Thle3}.
\end{enumerate}

\paragraph{Extreme values on a dependency graph}
By a dependency graph, we mean a graph $G=(V,E)$ and a collection of real random variables $X_{\mathbf{i}}, \mathbf{i}\in V$ (not necessarily stationary)  which satisfy the following property: for any pair of disjoint sets $A_1,A_2\subset V$ such that no edge in $E$ has one endpoint in $A_1$ and the other in $A_2$, the $\sigma$-field $\sigma(X_{\mathbf{i}}, \mathbf{i}\in A_1)$ and $\sigma(X_{\mathbf{i}}, \mathbf{i}\in A_2)$ are mutually independent. Introduced by Petrovskaya and Leontovitch in \cite{PL}, this concept was applied by Baldi and Rinott (e.g. \cite{BR}) to obtain central limit theorems and normal approximations. Furthermore, Arratia \textit{et al.} give a Poisson approximation of a sum of (non independent) Bernoulli random variables for a random field (see Theorem 1 in \cite{AGG}). We write their result in our context to approximate the number of exceedances on a dependency graph by a Poisson random variable.

First, we give some notations. We denote by $|V|$ the number of vertices of $G=(V,E)$, $D$ the maximal degree and $J\subset\RR$ a finite union of disjoint intervals. Let $\mathbf{U}'_V(J)$ be the number of exceedances i.e. \[\mathbf{U}'_V(J)=\sum_{\mathbf{i}\in V}\mathbb{1}_{X_{\mathbf{i}}\in J}\] and $p_{\mathbf{i}}=\PP(X_{\mathbf{i}}\in J)$, $p_{\mathbf{ij}}=\PP(X_{\mathbf{i}}\in J, X_{\mathbf{j}}\in J)$ for all $\mathbf{i}\in V$ and $\mathbf{j}\in V(\mathbf{i})-\{\mathbf{i}\}$ where $V(\mathbf{i})$ is the set of neighbors of $\mathbf{i}$ i.e. \begin{equation}\label{defVgraph}V(\mathbf{i}) = \{\mathbf{j}\in V, (\mathbf{i},\mathbf{j})\in E\}\cup\{\mathbf{i}\}.\end{equation} Let us consider a Poisson random variable $Z$ of mean \[\boldsymbol{\mu}_J = \EE[Z]=\EE[\mathbf{U}'_V(J)] = \sum_{\mathbf{i}\in V}\PPP{X_{\mathbf{i}}\in J}.\] Chen-Stein method can be applied to approximate the number of occurrences of dependent events by a Poisson random variable (e.g. \cite{AGG}). In particular, this is a powerful tool to derive some results in extreme value theory for a sequence of real random variables (e.g. \cite{S}). We write below a slightly modified version of Theorem 1 of \cite{AGG} to derive an upper bound of the total variation distance between the number of exceedances $\mathbf{U}'_V(J)$ and its Poisson approximation $Z$ for a dependency graph.

\begin{Prop}(Arratia \textit{et al.} 1989)
\label{steingraph}
Let $p(V) = \sup_{i\in V}p_{\mathbf{i}}$ and $q(V)^2 = \sup_{(\mathbf{i},\mathbf{j})\in E} p_{{\mathbf
i}\mathbf{j}}$. Then \begin{equation}\label{totalvar}\sup_{A\subset \NN}\left|\PP(\mathbf{U}'_V(J)\in A)-\PP(Z\in A) \right|\leq 2D\cdot |V|\cdot\left(p(V)^2+q(V)^2\right).\end{equation} In particular, for all $r\in\NN^*$, we get
\begin{equation}\label{totalvar0}\left|\PP(\mathbf{U}'_V(J)\leq r-1)-e^{-\boldsymbol{\mu}_J}\sum_{k=0}^{r-1}\frac{\boldsymbol{\mu}_J^k}{k!} \right|\leq 2D\cdot |V|\cdot\left(p(V)^2+q(V)^2\right).\end{equation}

\end{Prop}

\begin{prooft}{Proposition \ref{steingraph}}
The upper bound \eqref{totalvar0} is a direct consequence of \eqref{totalvar}. From Theorem 1 of \cite{AGG}, we get \begin{equation}\label{majstein}\sup_{A\subset \NN}\left|\PP(\mathbf{U}'_V(J)\in A)-\PP(Z\in A) \right|\leq 2(b_1+b_2+b_3)\end{equation} where \[b_1=\sum_{\mathbf{i}\in V}\sum_{\mathbf{j}\in V(\mathbf{i})}p_{\mathbf{i}}p_{\mathbf{j}},\hspace{0.3cm} b_2=\sum_{\mathbf{i}\in V}\sum_{\mathbf{i}\neq \mathbf{j}\in V(\mathbf{i})}p_{\mathbf{ij}} \text{ and } b_3=\sum_{i\in V}\EE\left[\EE\left[X_{\mathbf{i}}-p_{\mathbf{i}}|\sigma(X_{\mathbf{j}}:\mathbf{j}\not\in V(\mathbf{i})) \right] \right].\] Since $|V(\mathbf{i})|\leq D+1$, we obtain $b_1\leq |V|\cdot D\cdot p(V)^2$ and $b_2\leq |V|\cdot D\cdot q(V)^2$. Moreover, using the fact that if $\mathbf{j}\not\in V(\mathbf{i})$, the random variable $X_{\mathbf{j}}$ is independent of $X_{\mathbf{i}}$, we get $b_3=0$. We then deduce \eqref{totalvar} from \eqref{majstein}. 
\end{prooft}
Central limit theorems in geometric probability have been deduced from normal approximation on a dependency graph by a  
discretization technique (see e.g. \cite{AB}). In the same spirit, we derive Lemma \ref{Thle2} from Proposition \ref{steingraph}. We need first to explain how we construct the dependency graph from our random tessellation.

\paragraph{Construction of the dependency graph}
 We define a graph $G_{\rho} = (V_{\rho}, E_{\rho})$ as follows. The set $V_{\rho}$ consists of the sub-cubes $\mathbf{i}$ ($|V_{\rho}| = N_{\rho}$) which cover $\mathbf{W}_{\rho}$ whereas an edge $(\mathbf{i}, \mathbf{j})\in E_{\rho}$ if $d(\mathbf{i}, \mathbf{j})\leq R$ where $R$ is introduced in \textsc{Condition 1}. The maximal degree $D_{\rho}$ of this graph satisfies \begin{equation}\label{majDegree}D_{\rho} \leq (2R+1)^{d}.\end{equation} For all $\mathbf{i}\in V_{\rho}$, we define the random variable $X_{\mathbf{i}}$ as 
\begin{equation}\label{defXi} X_{\mathbf{i}} = M_{f,\mathbf{i}}.\end{equation}
From \textsc{Condition 1}, conditional on $A_{\rho}$, the graph $G_{\rho}$ and the collection $(M_{f, \mathbf{i}})_{\mathbf{i}\in V_{\rho}}$ define a dependency graph.

\begin{prooft}{Lemma \ref{Thle2}}
We apply Proposition \ref{steingraph} to $X_{\mathbf{i}} = M_{f, \mathbf{i}}$ and $J=(v_{\rho},\infty)$. It is enough to derive upper bounds of $\PP(M_{f, \mathbf{i}}>v_{\rho}|A_{\rho})$ and $\PP(M_{f, \mathbf{i}}>v_{\rho},M_{f, \mathbf{j}}>v_{\rho}|A_{\rho})$. According to \eqref{defXi}, we get 
\[\PP(M_{f, \mathbf{i}}>v_{\rho}) = \PPP{\bigcup_{\underset{z(C)\in \mathbf{i}\cap \mathbf{W}_{\rho}}{C\in\mathfrak{m},}}\{f(C)>v_{\rho}\}}  \leq \EEE{\sum_{\underset{z(C)\in \mathbf{i}\cap \mathbf{W}_{\rho}}{C\in\mathfrak{m},}}\mathbb{1}_{f(C)>v_{\rho}}}.\] Since $f$ is translation invariant and $\lambda_d(\mathbf{i})=\lambda_d(W)\cdot\rho/N_{\rho}$, we deduce from \eqref{campbell} that \begin{equation}\label{majcubei}\PP(M_{f, \mathbf{i}}>v_{\rho})\leq \frac{1}{N_{\rho}}\cdot\lambda_d(W)\cdot\rho\cdot\PP(f(\cell)>v_{\rho}).\end{equation} Using the trivial inequalities  $\PP(M_{f, \mathbf{i}}>v_{\rho}|A_{\rho})\leq \PP(M_{f, \mathbf{i}}>v_{\rho})/\PP(A_{\rho})$ and  $\lambda_d(W)\cdot\rho\cdot\PP(f(\cell)>v_{\rho})\leq G_1(\rho)+\tau$ where $G_1(\rho)$ is defined in \eqref{defG1general}, we obtain \begin{equation}\label{majorder1}
p_{\mathbf{i}}:=\PP(M_{f, \mathbf{i}}>v_{\rho}|A_{\rho}) \leq \frac{G_1(\rho) + \tau }{\PP(A_{\rho})N_{\rho}}.
\end{equation}  

Moreover, for any $\mathbf{i}\in V_{\rho}$ and $\mathbf{j}\in V_{\rho}(\mathbf{i})-\{\mathbf{i}\}$, we get
\begin{multline}\label{majcouple}\PP(M_{f, \mathbf{i}}>v_{\rho}, M_{f, \mathbf{j}}>v_{\rho}) = \PPP{\bigcup_{\underset{ z(C_1)\in \mathbf{i}\cap \mathbf{W}_{\rho}}{C_1\in\mm,}}\bigcup_{\underset{z(C_2)\in \mathbf{j}\cap \mathbf{W}_{\rho}}{C_2\in\mm,}}\{f(C_1)>v_{\rho}, f(C_2)>v_{\rho}\} }\\ \leq \EEE{\sum_{\underset{z(C_1),z(C_2)\in V_{\rho}(\mathbf{i})}{(C_1,C_2)_{\neq}\in\mm^2}}\mathbb{1}_{f(C_1)>v_{\rho}, f(C_2)>v_{\rho}}}\end{multline} where $(C_1,C_2)_{\neq}\in\mm^2$ means that $(C_1,C_2)$ is a couple of distinct cells. With the slight abuse of notation, we will write in the rest of the paper $V_{\rho}(\mathbf{i})$ for the union of the sub-cubes $\bigcup_{\mathbf{j}\in V_{\rho}(\mathbf{i})}\mathbf{j}$. 

Besides, the set of neighbors $V_{\rho}(\mathbf{i})$ can be re-written as $V_{\rho}(\mathbf{i}) = \{\mathbf{j}\in V_{\rho}, d(\mathbf{i}, \mathbf{j})\leq R\}$. Hence  $V_{\rho}(\mathbf{i})$ is a convex union of disjoint sub-cubes of volume $\lambda_d(W)\cdot\rho/N_{\rho}$, which are at most $(2R+1)^d$, and can be included in the cube $\mathfrak{C}_{\rho}$ defined in \eqref{defcubeC} up to a translation. Since $f$ is translation invariant, we obtain \begin{equation}
\label{majneighbor}
\EEE{\sum_{\underset{z(C_1),z(C_2)\in V_{\rho}(\mathbf{i})}{(C_1,C_2)_{\neq}\in\mm^2}}\mathbb{1}_{f(C_1)>v_{\rho}, f(C_2)>v_{\rho}}} \leq \frac{G_2(\rho)}{N_{\rho}}.
\end{equation}
Using the fact that $\PP(M_{f, \mathbf{i}}>v_{\rho}, M_{f, \mathbf{j}}>v_{\rho}|A_{\rho})\leq \PP(M_{f, \mathbf{i}}>v_{\rho}, M_{f, \mathbf{j}}>v_{\rho})/\PP(A_{\rho})$  we deduce from \eqref{majcouple} that
\begin{equation}
\label{majorder2}
p_{\mathbf{i}\mathbf{j}}:=\PP(M_{f, \mathbf{i}}>v_{\rho}, M_{f, \mathbf{j}}>v_{\rho}|A_{\rho}) \leq \frac{G_2(\rho)}{\PP(A_{\rho})N_{\rho}}.
\end{equation}

From \eqref{totalvar0} written for the conditional probability $\cdot|A_{\rho}$, \eqref{majDegree}, \eqref{majorder1}, \eqref{majorder2} and the fact that $|V_{\rho}| = N_{\rho}$, we get
\[\left|\PP(U'_{V_{\rho}}(v_{\rho})\leq r-1|A_{\rho}) - e^{-\mu_{\rho}}\sum_{k=0}^{r-1}\frac{\mu_{\rho}^k}{k!} \right| \leq \frac{2(2R+1)^d}{\PP(A_{\rho})^2}\cdot \left(\frac{(G_1(\rho)+\tau)^2}{N_{\rho}} + \PP(A_{\rho})G_2(\rho)\right). \] The rate of convergence \eqref{majUVgamma} results directly from the previous upper bound and the fact that $\PP(A_{\rho})$ and $G_1(\rho)$ converge respectively to 1 and 0 according to \textsc{Condition (TCP)} and \textsc{Condition 1}. 
\end{prooft}

We prove below Lemmas \ref{Thle1} and \ref{Thle3}. 

\begin{prooft}{Lemma \ref{Thle1}}
Let us notice that Lemma \ref{Thle1} is trivial when $r=1$. More generally, for all $r\in\NN^*$, we have
\begin{equation}
\label{diffproba}
\left|\PPP{U_{\rho}(v_{\rho})\leq r-1} - \PPP{U'_{V_{\rho}}(v_{\rho})\leq r-1} \right| \leq 2\PPP{U_{\rho}(v_{\rho})\neq U'_{V_{\rho}}(v_{\rho})}.
\end{equation}
According to \eqref{defUgammaUV},  the above random variables differ if and only if there are at least two exceedances in the same sub-cube $\mathbf{i}$ i.e.
\begin{multline}
\label{MajLemma1}
\PPP{U_{\rho}(v_{\rho})\neq U'_{V_{\rho}}(v_{\rho})} = \PPP{\bigcup_{\mathbf{i}\in V_{\rho}}\bigcup_{\underset{z(C_1),z(C_2)\in \mathbf{i}\cap \mathbf{W}_{\rho}}{(C_1,C_2)_{\neq}\in\mm^2,}}\{f(C_1)>v_{\rho}, f(C_2)>v_{\rho}\}}\\\leq \sum_{\mathbf{i}\in V_{\rho}}\EEE{\sum_{\underset{z(C_1),z(C_2)\in\mathbf{i}\cap \mathbf{W}_{\rho}}{(C_1,C_2)_{\neq}\in\mm^2}}\mathbb{1}_{f(C_1)>v_{\rho}, f(C_2)>v_{\rho}}}.
\end{multline}
 Since $|V_{\rho}| = N_{\rho}$, the right-hand side is bounded by $G_2(\rho)$ thanks to \eqref{majneighbor}. This shows that  $\PPP{U_{\rho}(v_{\rho})\neq U'_{V_{\rho}}(v_{\rho})}\leq G_2(\rho)$ and consequently we deduce \eqref{majThle1} from \eqref{diffproba}.
\end{prooft}

\begin{prooft}{Lemma \ref{Thle3}}
From \eqref{deftaugamma} and the triangle inequality, we get
\begin{equation}\label{triangleinequalitymu}|\mu_{\rho}-\tau| \leq |\EEE{U_{\rho}(v_{\rho})}-\tau| + \EEE{U_{\rho}(v_{\rho})-U'_{V_{\rho}}(v_{\rho})}\end{equation} where $U_{\rho}(v_{\rho})\geq U'_{V_{\rho}}(v_{\rho})$ a.s. According to \eqref{meanexceedance} and \eqref{defG1general}, we obtain that \begin{equation}\label{majU'1}\left|\EE[U_{\rho}(v_{\rho})]-\tau \right| =  G_1(\rho). \end{equation} 
To give an upper bound of the second term of the right-hand side of \eqref{triangleinequalitymu}, we use the fact that the family $V_{\rho}$ covers $\mathbf{W}_{\rho}$. Intuitively, the number of exceedance sub-cubes $U_{{V}_{\rho}}$ can be approximated by the number of exceedance cells $U_{\rho}(v_{\rho})$ when  $G_2(\rho)$ is negligible. We justify this fact below. From \eqref{defUgammaUV}, we obtain a.s. that \begin{multline}\label{diffexceedance}U_{\rho}(v_{\rho})-U'_{V_{\rho}}(v_{\rho}) = \sum_{\mathbf{i}\in V_{\rho}}\sum_{\underset{z(C)\in\mathbf{i}\cap \mathbf{W}_{\rho}}{C\in\mm,}}\mathbb{1}_{f(C)>v_{\rho}}-\mathbb{1}_{M_{f,\mathbf{i}}>v_{\rho}} = \sum_{\mathbf{i}\in V_{\rho}}\left(\sum_{\underset{z(C)\in\mathbf{i}\cap \mathbf{W}_{\rho}}{C\in\mm,}}\mathbb{1}_{f(C)>v_{\rho}}-1\right)\mathbb{1}_{M_{f,\mathbf{i}}>v_{\rho}}\\
\leq \sum_{\mathbf{i}\in V_{\rho}}\sum_{\underset{z(C_1),z(C_2)\in\mathbf{i}\cap \mathbf{W}_{\rho}}{(C_1,C_2)_{\neq}\in\mm^2}}\mathbb{1}_{f(C_1)>v_{\rho}, f(C_2)>v_{\rho}} .\end{multline} The last inequality comes from the fact that if there is 0 or 1 exceedance cell, the sums inside the expectations are null. Otherwise, if the number of exceedances is $k\geq 2$, we use that fact that $k-1\leq \frac{k(k-1)}{2}=:\binom{k}{2}$ which is the number of exceedance couples.

 Taking the means in \eqref{diffexceedance} and using the fact that the mean of the right-hand side of \eqref{diffexceedance} is bounded by $G_2(\rho)$ as in the proof of  Lemma \ref{Thle1}, we get
  \begin{equation}
\label{majU'2}
\EEE{U_{\rho}(v_{\rho})-U'_{V_{\rho}}(v_{\rho})} \leq G_2(\rho)
\end{equation}
From \eqref{triangleinequalitymu}, \eqref{majU'1} and \eqref{majU'2} we obtain that  $|\mu_{\rho} - \tau|$ is lower than $G_1(\rho)+ G_2(\rho)$. 
\end{prooft}

\subsection{Proof of Theorem \ref{PPTh}}

By Kallenberg's theorem (see Proposition 3.22, p. 156 in \cite{R}, see also the proof of Theorem 2.1.2 in \cite{HF}) it is enough to check that:
\begin{itemize}
\item For all Borel subset $B\subset W$ and ${}_*x<s\leq t\leq x^*$
\begin{equation}\label{PPpart1}\EEE{\#\Phi_{\rho}\cap (B\times (s,t])} \conv[\rho]{\infty}\EEE{\#\Phi\cap (B\times (s,t]) } \end{equation}
\item For all $\mathscr{P}=\bigcup_{l=1}^LB^{(l)}\times (s_l,t_l]$ where $B^{(l)}$ is the intersection of $W$ and a rectangular solid in $\mathbf{C}^{(W)}$ and ${}_*x<s_l\leq t_l\leq x^*$
 \begin{equation}\label{PPpart2}\PPP{\#\Phi_{\rho}\cap \mathscr{P} = 0} \conv[\rho]{\infty}\PPP{\#\Phi\cap \mathscr{P} = 0}\end{equation}
\end{itemize}

\begin{prooft}{\eqref{PPpart1}} From \eqref{campbell}, we have
\begin{equation*} \EEE{\#\Phi_{\rho}\cap (B\times (s,t])}   =  \EEE{\sum_{\underset{z(C)\in \mathbf{B}_{\rho}}{C\in\mm,}}\mathbb{1}_{a_{\rho}s + b_{\rho}<f(C)\leq a_{\rho}t + b_{\rho}} }
 = \lambda_d(\mathbf{B}_{\rho})\cdot \left( \PPP{f(\cell)>a_{\rho}s +b_{\rho}} -  \PPP{f(\cell)>a_{\rho}t +b_{\rho}} \right)\end{equation*} where we recall that $\mathbf{B}_{\rho} = \rho^{1/d}B$. According to the trivial equality $\lambda_d(\mathbf{B}_{\rho}) = \frac{\lambda_d(B)}{\lambda_d(W)}\cdot \lambda_d(W)\cdot \rho$ and the fact that $\lambda_d(\mathbf{W}_{\rho})\cdot \PPP{f(\cell)>v_{\rho}(t)}$ converges to $\tau(t)$ for all $t\in\RR$, we get \begin{equation}\label{convexpPP}\EEE{\#\Phi_{\rho}\cap (B\times (s,t])}\conv[\rho]{\infty}\frac{\lambda_d(B)}{\lambda_d(W)}\cdot (\tau(s)-\tau(t)) = \EEE{\#\Phi\cap (B\times (s,t])}\end{equation}  and consequently we obtain \eqref{PPpart1}.
 \end{prooft}
 
\begin{prooft}{\eqref{PPpart2}} We can write $\mathscr{P}$ as a disjoint union of strips i.e.
\begin{equation}
\label{notP}
\mathscr{P}=\bigsqcup_{l=1}^LB^{(l)}\times J^{(l)}
\end{equation} such that the Borel subsets $B^{(l)}\subset W$ are disjoint and  such that $J^{(l)}$ is a finite union of half-open intervals for all $l=1,\ldots, L$. The following lemma shows that it is enough to investigate the case where $\mathscr{P}$ is a strip.

\begin{Le}
\label{lemmadisjoint}
Let $\mathscr{P}$ be as in \eqref{notP}. With the same hypothesis as in Theorem \ref{PPTh}, we have
\begin{equation}\label{convdisjoint}\PPP{\#\Phi_{\rho}\cap \mathscr{P} = 0} - \prod_{l=1}^L\PPP{\#\Phi_{\rho}\cap (B^{(l)}\times J^{(l)}) = 0} \conv[\rho]{\infty}0.\end{equation}
\end{Le}
The proof of Lemma \ref{lemmadisjoint} is postponed at the end of the subsection. Thanks to Lemma \ref{lemmadisjoint}, we can assume that $\mathscr{P}$, defined in \eqref{notP}, is only a strip i.e. $\mathscr{P}=B\times J$ where $J$ is a finite union of half-open intervals. Without loss of generality, we can assume that these intervals are disjoint i.e. \begin{equation}\label{defJ}J=\bigsqcup_{j=1}^k (s_j,t_j]\end{equation} with ${}_*x<s_j\leq t_j\leq x^*$ and $s_j\leq t_{j+1}$, $j=1,\ldots, k$. In the same spirit as in the proof of Theorem \ref{realtyp}, we introduce two random variables that are
\begin{equation} \label{defUPP}\mathscr{U}_{\rho}(B\times J) = \#\Phi_{\rho}\cap (B\times J) = \sum_{\underset{z(C)\in \mathbf{B}_{\rho}}{C\in\mm,}}\mathbb{1}_{a_{\rho}^{-1}(f(C) - b_{\rho} )\in J} \text{ and } \mathscr{U}'_{V_{\rho}}(B\times J) = \sum_{\mathbf{i}\in V_{\rho}}\mathbb{1}_{a_{\rho}^{-1}(M_{f,\mathbf{i}}(B) - b_{\rho} )\in J}  \end{equation} where \[M_{f,\mathbf{i}}(B) = \max_{\underset{z(C)\in\mathbf{i}\cap \mathbf{B}_{\rho}}{C\in\mm,}}f(C).\] In particular, $\mathscr{U}_{\rho}(W\times (s, \infty)) = U_{\rho}(v_{\rho}(s))$  and $\mathscr{U}'_{V_{\rho}}(W\times (s, \infty)) = U'_{V_{\rho}}(v_{\rho}(s))$ where $U_{\rho}(v_{\rho}(s))$ and $ U'_{V_{\rho}}(v_{\rho}(s))$ have been defined in \eqref{defUgammaUV}. We denote by $\boldsymbol{\mu}_{\rho}(B\times J)$ the mean of $\mathscr{U}'_{V_{\rho}}(B\times J)$ i.e. 
\[\boldsymbol{\mu}_{\rho}(B\times J) = \EEE{\mathscr{U}'_{V_{\rho}}(B\times J)} = \sum_{\mathbf{i}\in V_{\rho}}\PPP{a_{\rho}^{-1}(M_{f,\mathbf{i}}(B)-b_{\rho})\in J}.\] As in the proof of Theorem \ref{realtyp}, we subdivide the proof into three steps. More precisely, we show that
\begin{subequations}
\begin{equation}
\label{PPpart2.1}
\PPP{\mathscr{U}_{\rho}(B\times J) = 0} - \PPP{\mathscr{U}'_{V_{\rho}}(B\times J) = 0} \conv[\rho]{\infty}0
\end{equation}
\begin{equation}
\label{PPpart2.2}
\PPP{\mathscr{U}'_{V_{\rho}}(B\times J) = 0} - e^{-\boldsymbol{\mu}_{\rho}(B\times J)}\conv[\rho]{\infty}0
\end{equation}
\begin{equation}
\label{PPpart2.3}
\boldsymbol{\mu}_{\rho}(B\times J) \conv[\rho]{\infty} \nu(B\times J). 
\end{equation}
\end{subequations}
Let us notice that the convergences \eqref{PPpart2.1},\eqref{PPpart2.2} and \eqref{PPpart2.3} are generalisations of Lemmas \ref{Thle1}, \ref{Thle2} and \ref{Thle3} respectively. For the proof of \eqref{PPpart2.1}, it is enough to show that $\PPP{\mathscr{U}_{\rho}(B\times J)\neq \mathscr{U}'_{V_{\rho}}(B\times J)}$ converges to 0 as $\rho$ goes to infinity. Since $\mathscr{U}_{\rho}(B\times J)\geq \mathscr{U}'_{V_{\rho}}(B\times J)$ for all Borel subsets, we have
\begin{multline*}
\label{majdiffPP}
\PPP{\mathscr{U}_{\rho}(B\times J)\neq \mathscr{U}'_{V_{\rho}}(B\times J)} \leq \sum_{j=1}^k\PPP{\mathscr{U}_{\rho}(B\times (s_j,t_j])\neq \mathscr{U}'_{V_{\rho}}(B\times (s_j,t_j])}\\
\leq \sum_{j=1}^k\PPP{\mathscr{U}_{\rho}(W\times (s_j,\infty))\neq \mathscr{U}'_{V_{\rho}}(W\times (s_j,\infty))} = \sum_{j=1}^k\PPP{U_{\rho}(v_{\rho}(s_j))\neq U'_{V_{\rho}(v_{\rho}(s_j))}}. 
\end{multline*} 
Bounding  as in \eqref{MajLemma1} and proceeding along the same lines as in the proof of Lemma \ref{Thle1}, we show that the right-hand side converges to 0.

Secondly, we prove \eqref{PPpart2.2}. In the same spirit as in the proof of Lemma \ref{Thle2}, we apply Proposition \ref{steingraph} conditional on $A_{\rho}$  to $X_{\mathbf{i}} = a_{\rho}^{-1}(M_{f,\mathbf{i}}(B)-b_{\rho})$ and $J=\bigsqcup_{j=1}^k (s_j,t_j]$. Let $\mathbf{i}\in V_{\rho}$ and $\mathbf{j}\in V_{\rho}(\mathbf{i})-\{\mathbf{i}\}$. Using the fact that $M_{f,\mathbf{i}}(B)\leq M_{f,\mathbf{i}}$, we get \[p_{\mathbf{i}} = \PPP{ a_{\rho}^{-1}(M_{f,\mathbf{i}}(B)-b_{\rho})\in J|A_{\rho}} \leq \PPP{M_{f,\mathbf{i}}(B)>v_{\rho}(s_1)|A_{\rho}} = O(N_{\rho}^{-1})\] according to \eqref{majcubei}. Moreover
\begin{multline*} p_{\mathbf{ij}} = \PPP{a_{\rho}^{-1}(M_{f,\mathbf{i}}(B)-b_{\rho})\in J, a_{\rho}^{-1}(M_{f,\mathbf{j}}(B)-b_{\rho})\in J|A_{\rho}}\leq \PPP{M_{f,\mathbf{i}}(B)>v_{\rho}(s_1), M_{f,\mathbf{i}}(B)>v_{\rho}(s_1)|A_{\rho}}\\ = O\left(G_2(\rho)\cdot N_{\rho}^{-1} \right)\end{multline*} according to \eqref{majorder2}. We deduce \eqref{PPpart2.2} from the previous inequalities and Proposition \ref{steingraph}. 

Finally, we prove \eqref{PPpart2.3}. According to  \eqref{defJ} and \eqref{defUPP}, we have a.s.
\[\mathscr{U}_{\rho}(B\times J) = \sum_{j=1}^k \#\Phi_{\rho}\cap (B\times (s_j,t_j]).\] Taking the expectations in the previous equality, we deduce from \eqref{convexpPP} that \begin{equation}\label{PPpart2.3.1}  \EEE{\mathscr{U}_{\rho}(B\times J)} \conv[\rho]{\infty}\frac{\lambda_d(B)}{\lambda_d(W)}\sum_{j=1}^k(\tau(s_j)-\tau(t_j))=\nu(B\times J).\end{equation} Moreover
\begin{equation}\label{PPpart2.3.2}\EEE{\mathscr{U}_{\rho}(B\times J)} - \boldsymbol{\mu}_{\rho}(B\times J) = \EEE{\mathscr{U}_{\rho}(B\times J) - \mathscr{U}'_{V_{\rho}}(B\times J)} \leq \EEE{U_{\rho}(v_{\rho}(s_1)) - U'_{V_{\rho}}(v_{\rho}(s_1))}\end{equation} converges to 0 according to \eqref{majU'2}. We deduce \eqref{PPpart2.3} from \eqref{PPpart2.3.1} and \eqref{PPpart2.3.2}. 

\subparagraph{Conclusion of the proof of \eqref{PPpart2}.} According to \eqref{PPpart2.1}, \eqref{PPpart2.2},   \eqref{PPpart2.3} and the fact that $\mathscr{U}_{\rho}(B\times J) = \#\Phi_{\rho}\cap(B\times J)$, we deduce that \[\PPP{\#\Phi_{\rho}\cap (B\times J) = 0} \conv[\rho]{\infty} e^{-\nu(N\times J)} =  \PPP{\#\Phi\cap (B\times J) = 0}\] and consequently we obtain \eqref{PPpart2}. 
\end{prooft}

The end of the subsection is devoted to the proof of Lemma \ref{lemmadisjoint}.

\begin{prooft}{Lemma \ref{lemmadisjoint}}
Let $\mathscr{P}=\bigsqcup_{l=1}^LB^{(l)}\times J^{(l)}$ and  $B^{(l)}= B_l\cap W$  such that the rectangular solids $B_l\subset \mathbf{C}^{(W)}$ are disjoint. First, we introduce some notations. We denote by $V_{\rho}(B^{(l)})$, $S_{\rho}(B^{(l)})$ and  $V^{\circ}_{\rho}(B^{(l)})$ respectively  the sets \[\left\{ \begin{split} & V_{\rho}(B^{(l)}) = \{\mathbf{i}\in V_{\rho}, \mathbf{i}\cap B_l\neq\varnothing\}\\
& S_{\rho}(B^{(l)}) = \{\mathbf{i}\in V_{\rho}, \mathbf{i}\cap\partial B_l\neq\varnothing\} \\
&  V^{\circ}_{\rho}(B^{(l)}) = \{\mathbf{i}\in V_{\rho}(B^{(l)}), d(\mathbf{i}, S_{\rho}(B^{(l)}))>R \} \end{split}\right..\]  Finally, we denote by $\mathscr{U}'_{V^{\circ}_{\rho}}(B^{(l)}\times J^{(l)})\leq \mathscr{U}_{\rho}(B^{(l)}\times J^{(l)})$ the number of exceedances in  $V^{\circ}_{\rho}(B^{(l)})$ i.e. 
\[\mathscr{U}'_{V^{\circ}_{\rho}}(B^{(l)}\times J^{(l)}) = \sum_{\mathbf{i}\in  V^{\circ}_{\rho}(B^{(l)})}\mathbb{1}_{a_{\rho}^{-1}(M_{f,\mathbf{i}}(B^{(l)})-b_{\rho})\in J^{(l)}}. \]  
Let $l\in \{1,\ldots, L\}$ be fixed. Since $B_l$ is a rectangular solid in $\mathbf{C}^{(W)}$ which is covered with at most $N_{\rho}$ sub-cubes $\mathbf{i}$, we have $\# S_{\rho}(B^{(l)}) \leq c\cdot N_{\rho}^{(d-1)/d}$. This shows that
\[\PPP{\mathscr{U}'_{V^{\circ}_{\rho}}(B^{(l)}\times J^{(l)}) \neq \mathscr{U}'_{V_{\rho}}(B^{(l)}\times J^{(l)})} \leq \# S_{\rho}(B^{(l)})\cdot \PPP{M_{f,\mathbf{i}}>v_{\rho}} = O\left( N_{\rho}^{-1/d} \right) \] according to \eqref{majcubei} and \textsc{Condition (TCP)}. Thanks to \eqref{PPpart2.1}, we deduce that
 \begin{equation}\label{approxindependence}\PPP{\mathscr{U}_{\rho}(B^{(l)}\times J^{(l)})=0} - \PPP{\mathscr{U}'_{V^{\circ}_{\rho}}(B^{(l)}\times J^{(l)})=0}  \conv[\rho]{\infty}0.\end{equation}

Moreover, conditional on $A_{\rho}$, the random variables $\mathscr{U}'_{V^{\circ}_{\rho}}(B^{(l)}\times J^{(l)})$, $l=1,\ldots, L$ are independent since the rectangular solids $B_l$, $l=1,\ldots, L$ are at distance higher than $R$. In  particular, we get
\[\PPP{\bigcap_{l=1}^L\left\{\left.\mathscr{U}'_{V^{\circ}_{\rho}}(B^{(l)}\times J^{(l)}) = 0\right\}\right|A_{\rho}} = \prod_{l=1}^L\PPP{\left.\mathscr{U}'_{V^{\circ}_{\rho}}(B^{(l)}\times J^{(l)}) = 0\right|A_{\rho}}.\] Lemma \ref{lemmadisjoint} is a consequence of the previous equality, the convergence \eqref{approxindependence} and the fact that \[\PPP{\#\Phi_{\rho}\cap \mathscr{P} = 0} = \PPP{\bigcap_{l=1}^L\{\mathscr{U}_{\rho}(B^{(l)}\times J^{(l))}) = 0\}}.\]
\end{prooft}

\begin{Rk}
\label{Rkdisjoint}
When \textsc{Condition 2} does not hold, Lemma \ref{lemmadisjoint} remains true when $\mathscr{P}=\bigsqcup_{l=1}^LB^{(l)}\times (s_l,\infty)$. This comes from the fact that the left-hand side of \eqref{PPpart2.1} equals 0 when $J=(s,\infty)$. In the same spirit, we can show that
if $B^{(1)}, \ldots, B^{(l)}$, $1\leq l\leq L$ is a set of $L\geq 1$ disjoint Borel subsets included in $W$, we have:
\begin{equation}\label{Rkindependence}\PPP{M_{f, \mathbf{W}_{\rho}} \leq v_{\rho}} - \prod_{l=1}^L\PPP{M_{f,\mathbf{B}_{\rho}^{(l)}}\leq v_{\rho}} \conv[\rho]{\infty}0\end{equation} where $\mathbf{B}^{(l)}_{\rho}=\rho^{1/d}\mathbf{B}^{(l)}$, $1\leq l\leq L$. Let us note that the previous convergence holds for a threshold $v_{\rho}$ which is \textit{not necessarily} of the form $v_{\rho}=v_{\rho}(t)=a_{\rho}t+b_{\rho}$. We will use this remark in section \ref{sectionrealtyp2}. 
\end{Rk}

\begin{Rk}
The inequalities appearing in \eqref{defG1general}, \eqref{defG2general} and Theorem \ref{realtyp} have to be reversed when we deal with the $r$ smallest values. This fact will be extensively used in the rest of the paper. 
\end{Rk}

In the three following sections, we apply Theorem \ref{realtyp} to derive the asymptotic behaviours of the order statistics for different geometrical characteristics and random tessellations. For aesthetic reasons, we only investigate maxima and minima for the particular case $W=\mathbf{C}^{(W)}=[0,1]^d$ keeping in mind that these results can be generalized to order statistics and to any bounded set with $\lambda_d(W)\neq 0$. Up to a normalization, all the thresholds $v_{\rho}$ can be written as $v_{\rho}=v_{\rho}(t)=a_{\rho}t+b_{\rho}$ (excepted in section \ref{sectionGP}) so that Theorem \ref{PPTh} is also available.

\section{Extreme Values of a Poisson-Delaunay tessellation}
\label{sectionPDT}
Before applying Theorem \ref{realtyp} to different geometrical characteristics of a Poisson-Delaunay tessellation, we introduce some notations and preliminaries.

\paragraph{Notations}
\begin{itemize}
\item Let $z$ be a point in $\RR^d$ and $r$ be a positive real number. We denote by $B(z,r)$ and $S(z,r)$ the ball and the sphere of radius $r$ centered in $z$. When $z=0$ and $r=1$, we denote by $\mathbf{S}^{d-1}=S(0,1)$ the unit sphere. Moreover, we denote by $\kappa_d$ the volume of the unit ball i.e.
\[\kappa_d = \lambda_d(B(0,1)).\]
\item Let $C$ be a simplex in $\RR^d$. We denote respectively by $B(C)$, $S(C)$, $z(C)$ and $R(C)$ the circumball, the circumsphere, the circumcenter and the circumradius of $C$.
\item Let $k$ be an integer and $x_1,\ldots, x_k$ be $k$ points in $\RR^d$ and let $f:\mathcal{K}_d\rightarrow\RR$ be a measurable function.
\begin{itemize}
\item We denote by $\mathbf{x}_{1:k}$ the $k$-tuple $(x_1, \ldots, x_k)$ and by  $\{\mathbf{x}_{1:k}\}$ the set of points $\{x_1,\ldots, x_k\}$.
\item If $r$ is a positive number, we define $r\mathbf{x}_{1:k}=(rx_1, \ldots, rx_k)$ and $r\{\mathbf{x}_{1:k}\} = \{rx_1, \ldots, rx_k\}$. 
\item When $k=d+1$ and when the $d+1$ points $x_1,\ldots, x_{d+1}$ lie on a sphere, we denote by $\Delta(\mathbf{x}_{1:d+1})$ the convex hull of $x_1, \ldots, x_{d+1}$. Moreover, we define $f(\mathbf{x}_{1:d+1})$ as \[f(\mathbf{x}_{1:d+1}) = f\left(\Delta(\mathbf{x}_{1:d+1})\right).\] In particular, $B(\mathbf{x}_{1:d+1})$, $S(\mathbf{x}_{1:d+1})$, $z(\mathbf{x}_{1:d+1})$, $R(\mathbf{x}_{1:d+1})$ and $\lambda_d(\mathbf{x}_{1:d+1})$ are respectively the circumball, the circumsphere, the circumcenter, the circumradius and the volume of the simplex $\Delta(\mathbf{x}_{1:d+1})$. 
\item If $k\leq d+1$ and if $\{\mathbf{y}_{k+1:d+1}\}=\{y_{k+1}, \ldots, y_{d+1}\}$ is a set of $d+1-k$ points in $\RR^d$ such that $x_1,\ldots, x_{k}$ and $y_{k+1}, \ldots, y_{d+1}$ lie on a sphere, we denote by  $\Delta(\mathbf{x}_{1:k}, \mathbf{y}_{k+1:d+1})$ the convex hull of $x_1, \ldots, x_k, y_{k+1}, \ldots, y_{d+1}$. Moreover, we define $f(\mathbf{x}_{1:k}, \mathbf{y}_{k+1:d+1})$ as
\[f(\mathbf{x}_{1:k}, \mathbf{y}_{k+1:d+1}) = f\left(\Delta\left(\mathbf{x}_{1:k}, \mathbf{y}_{k+1:d+1}\right)\right).\]
\item Finally, we denote by $d\sigma(u)$ the uniform distribution over the unit sphere $\mathbf{S}^{d-1}$ and $d\sigma(\mathbf{u}_{1:d+1}) = d\sigma(u_1)\cdots d\sigma(u_{d+1})$.

\end{itemize}
\end{itemize}

\paragraph{Preliminaries}
Let $\chi$ be a locally finite subset of $\RR^d$ such that each subset of size $n<d+1$ are affinely independent and no $d+2$ points lie on a sphere. If $d+1$ points $x_1,\ldots, x_{d+1}$ of $\chi$ lie on a sphere that contains no point of $\chi$ in its interior, then the convex hull of $x_1,\ldots, x_{d+1}$ is called a cell. The set of such cells defines a  partition of $\RR^d$ into simplices and such partition is called the Delaunay tessellation. Such model is the key ingredient of the first algorithm for computing the minimum spanning tree \cite{ShHo}. It is extensively used in medical image segmentation \cite{SKSSS}, in finite element method to build meshes \cite{JGZ} and is a powerful tool for reconstructing a $3D$ set from a discrete point set \cite{Sch}. 

When $\chi = \mathbf{X}$ is a Poisson point process, we speak about Poisson-Delaunay tessellation and we denote this random tessellation by $\mm_{PDT}$. For each cell $C\in\mm_{PDT}$ which is a.s. a simplex, we define $z(C)$ as the circumcenter of $C$. The relation between the intensity $\gamma$ of $\mm_{PDT}$ and the intensity $\gamma_{\mathbf{X}}$ of the underlying Poisson point process is given by (see section 7 in \cite{M}) \[\gamma= \beta_d^{-1}\cdot \gamma_{\mathbf
X}\] where 
\begin{equation}
\label{defad} \beta_d = \frac{(d^3+d^2)\Gamma\left(\frac{d^2}{2} \right)\Gamma^d\left(\frac{d+1}{2} \right)}{\Gamma\left(\frac{d^2+1}{2} \right)\Gamma^{d}\left(\frac{d+2}{2} \right)2^{d+1}\pi^{\frac{d-1}{2}}}.   
\end{equation} To be in the framework of Theorem \ref{realtyp}, we assume (without loss of generality) that  $\gamma = 1$ i.e. 
\[\gamma_{\mathbf{X}} = \beta_d.\] Moreover, we partition the window $\mathbf{W}_{\rho}=\rho^{1/d}[0,1]^d$ into $N_{\rho}$ sub-cubes $\mathbf{i}\in V_{\rho}$ where we take
\[N_{\rho} = \left\lfloor\frac{\rho}{2\log\rho}\right\rfloor.\] To apply Theorem \ref{realtyp}, we first check \textsc{Condition 1} for any  measurable function $f:\mathcal{K}_d\rightarrow\RR$. To do it, we define the event $A_{\rho}$ (independent on $f$) as  \begin{equation}\label{defAgammadel}A_{\rho} = \bigcap_{i\in V_{\rho}}\{\mathbf{X}\cap \mathbf{i}\neq\varnothing\}.\end{equation}
\begin{Le}
\label{AssumptAgammadel}
Let $f:\mathcal{K}_d\rightarrow\RR$ be a measurable function. Then \textsc{Condition 1} is satisfied for $R=2\cdot\left(\lfloor \sqrt{d}\rfloor + 1\right)$ and for the event $A_{\rho}$ defined in \eqref{defAgammadel}.
\end{Le}

\begin{prooft}{Lemma \ref{AssumptAgammadel}}

We use the same arguments as in the proof of Proposition 3 in \cite{AB}. Let $\mathbf{i}\in V_{\rho}$ be a sub-cube in $\mathbf{W}_{\rho}$ and let $C\in \mm_{PDT}$ such that $z(C)\in\mathbf{i}$. Since a $d+1$-tuple of points of $\mathbf{X}$ is a Delaunay cell if and only if its circumball contains no point in its interior, we have $R(C) = \min_{x\in \mathbf{X}}\{|z(C)-x|\}$. Moreover, conditional on $A_{\rho}$, there exists a point $x_0$ in $\mathbf{X}\cap\mathbf{i}$. In particular, we have $|z(C)-x_0|\leq \sqrt{d}\cdot c_{\rho}$ where $c_{\rho}$ is the length of the sides of each sub-cube. Consequently, we obtain
\begin{equation} \label{eventcircumradius}R(C)\leq \sqrt{d}\cdot c_{\rho}.\end{equation}
This shows that the circumsphere $S(C)$ of $C$ is included in $V_{\rho}(\mathbf{i},D)$ where $D=\lfloor \sqrt{d}\rfloor +1$ and \[V_{\rho}(\mathbf{i}, D) = \{\mathbf{j}\in V_{\rho}, d(\mathbf{i},\mathbf{j})\leq D\}.\] Indeed if not, there exists a point $y\in S(C)$ such that $y$ is in a sub-cube $\mathbf{j}$ with $d(\mathbf{i},\mathbf{j})\geq D+1$. This shows that $|y-z(C)|> (\lfloor \sqrt{d}\rfloor +1)\cdot c_{\rho}$ and contradicts \eqref{eventcircumradius} since $R(C) = |y-z(C)|$. 

Since $S(C)$ is included in $V_{\rho}(\mathbf{i},D)$ for any cell $C\in\mm_{PDT}$ such that $z(C)\in \mathbf{i}$, this shows that $M_{f, \mathbf{i}}$ is $\sigma(\mathbf{X}\cap V_{\rho}(\mathbf{i},D))$ measurable. Because $d(A,B)>2D$ implies that $\{\mathbf{i}, d(\mathbf{i},A)<D\}$ and $\{\mathbf{i}, d(\mathbf{i},B)<D\}$ are disjoint and because $\mathbf{X}\cap \{\mathbf{i}, d(\mathbf{i},A)<D\}$ and $\mathbf{X}\cap \{\mathbf{i}, d(\mathbf{i},B)<D\}$ are independent, the $\sigma$-algebras $\sigma(M_{f,\mathbf{i}}, \mathbf{i}\in A)$ and $\sigma(M_{f,\mathbf{i}}, \mathbf{i}\in B)$ are independent, yielding $R=2D=2\cdot\left(\lfloor \sqrt{d}\rfloor + 1\right)$.

 Moreover the probability of the event $A_{\rho}$ converges to 1. Indeed, since $\mathbf{X}$ is a Poisson point process, we get
\begin{equation}\label{majPA}\PP(A_{\rho}^c) = \PPP{\bigcup_{\mathbf{i}\in V_{\rho}}\{\mathbf{X}\cap\mathbf{i}=\varnothing\}} \leq N_{\rho} e^{-\rho/N_{\rho}} = O\left((\log\rho)^{-1}\times\rho^{-1}\right).\end{equation}
\end{prooft}

Besides, the distribution function of the typical cell can be made explicit. Indeed, let $f:\mathcal{K}_d\rightarrow\RR$ be a translation invariant function on the set of convex bodies. An integral representation of $f(\cell)$, due to Miles \cite{Mi} (the proof can also be found in Theorem 10.4.4. of \cite{SW}), is given by 

\begin{equation}
\label{typicalcelldelaunay}
\EEE{f(\cell)} = \delta'_d\cdot\int_{0}^{\infty}\int_{(\mathbf{S}^{d-1})^{d+1}}r^{d^2-1}e^{-\delta_dr^d}\lambda_d(\mathbf{u}_{1:d+1})
f(r\mathbf{u}_{1:d+1})
d\sigma(\mathbf{u}_{1:d+1})dr
\end{equation} where 
\begin{equation}\label{defdelta}\delta'_d = (d+1)\cdot\beta_d \text{ and } \delta_d = \kappa_d\cdot\beta_{d}.\end{equation}
We recall that a $(d+1)$-tuple of points of $\mathbf{X}$ is a Delaunay cell if and only if its circumball contains no point of $\mathbf{X}$ in its interior. This justifies the exponential term since it is the probability that $\mathbf{X}\cap B(0,r)$ is empty. Thanks to \eqref{typicalcelldelaunay}, the typical cell can be built explicitly: it is a random simplex inscribed in the ball $B(0,r)$ such that the vector $\mathbf{u}_{1:d+1}$ is independent of $r$ and has a density proportional to the volume of the simplex $\Delta(\mathbf{u}_{1:d+1})$.

For practical reasons, we write below a generic lemma which gives an integral representation of the function $G_2(\cdot)$ defined in \eqref{defG2general}. To do it, we introduce some notations. As defined in \eqref{defG2general}, $G_2(\cdot)$ brings up two cells $\Delta_1, \Delta_2$ that are two different simplices such that $f(\Delta_i)>v_{\rho}$ and $z(\Delta_i)\in \mathfrak{C}_{\rho}$, $i=1,2$. The intersection of these cells is a $k$-dimensional simplex with $0\leq k\leq d-1$. Translating the circumcenter of the cell which has the largest circumradius say $\Delta_1$ at the origin, the cells can we written as $\Delta_1=\Delta(r\mathbf{u}_{1:d+1})$ and $\Delta_2=\Delta(r\mathbf{u}_{1:k}, \mathbf{y}_{k+1:d+1})$ with  $r\geq 0$, $u_1, \ldots, u_{d+1}\in \mathbf{S}^{d-1}$ and $y_{k+1}, \ldots, y_{d+1}\in\RR^d$. We  consider two properties $\mathscr{P}_1, \mathscr{P}_2$  that are 
\begin{subequations}
\begin{equation}
\label{delcond1}
\mathscr{P}_1: \vspace{0.3cm} f(r\mathbf{u}_{1:k}, \mathbf{y}_{k+1:d+1})>v_{\rho}, \vspace{0.3cm}   R(r\mathbf{u}_{1:k}, \mathbf{y}_{k+1:d+1})\leq r \text{ and } z(r\mathbf{u}_{1:k}, \mathbf{y}_{k+1:d+1})\in\CC.
\end{equation}
\begin{equation}
\label{delcond2}
\mathscr{P}_2: \vspace{0.3cm} y_j\not\in B(r\mathbf{u}_{1:d+1}) \text{ and } ru_j\not\in B(r\mathbf{u}_{1:k}, \mathbf{y}_{k+1:d+1}) \text{ for all } j=k+1,\ldots, d+1.
\end{equation}
\end{subequations} The first property concerns the cell $\Delta_2$ which has the smallest circumradius whereas the second property means that the two simplices are Delaunay cells. Moreover, we introduce the set
\begin{equation} \label{defE}E_{k,r,\mathbf{u}_{1:d+1}}  = \{\mathbf{y}_{k+1:d+1}\in (\RR^d)^{d+1-k} \text{ satisfying } \mathscr{P}_1 \text{ and } \mathscr{P}_2 \}.\end{equation} At last, in the same spirit as in \eqref{typicalcelldelaunay}, we consider the volume of the union of the two circumballs i.e. \begin{equation}\label{defunionvol}
\lambda^{(\cup)}_d(r, \mathbf{u}_{1:k}, \mathbf{y}_{k+1:d+1}) = \lambda_d\left( B(0,r)\cup B(r\mathbf{u}_{1:k}, \mathbf{y}_{k+1:d+1})\right).\end{equation} We are now prepared to state the generic lemma.

\begin{Le}
\label{genericdelaunay}
Let $\mm_{PDT}$ be a Poisson-Delaunay tessellation of intensity $\gamma = 1$. Then
\begin{equation}\label{majG2del}G_2(\rho) = 2\cdot\sum_{k=0}^d G_{2,k}(\rho)\end{equation} where 
\begin{equation}\label{defG2k} G_{2,k}(\rho) = \rho\int_0^{\infty}\int_{(\mathbf{S}^{d-1})^{d+1}}\int_{(\RR^d)^{d+1-k}} g_{2,k}(\rho, r, \mathbf{u}_{1:d+1},\mathbf{y}_{k+1:d+1})drd\sigma(\mathbf{u}_{1:d+1})d\mathbf{y}_{k+1:d+1}\end{equation} and \begin{equation}\label{defg2kdel}g_{2,k}(\rho, r, \mathbf{u}_{1:d+1},\mathbf{y}_{k+1:d+1}) = r^{d^2-1}  e^{-\beta_d \lambda^{(\cup)}_d(r, \mathbf{u}_{1:k}, \mathbf{y}_{k+1:d+1})}\lambda_d(\mathbf{u}_{1:d+1}) \mathbb{1}_{f(r\mathbf{u}_{1:d+1})>v_{\rho}}\mathbb{1}_{E_{k,r,\mathbf{u}_{1:d+1}}}(\mathbf{y}_{k+1:d+1}).\end{equation} 
\end{Le}

\begin{prooft}{Lemma \ref{genericdelaunay}}
This will be sketched since it in the same spirit as in the proof of \eqref{typicalcelldelaunay}. Considering that the intersection of the two Delaunay cells $\Delta_1$, $\Delta_2$ which appear in \eqref{defG2general} is a $k$-dimensional simplex with $0\leq k\leq d$ and assuming that $R(\Delta_1)\geq R(\Delta_2)$, we have
\begin{multline*}
\PPP{f(\cell)>v_{\rho}} = 2\sum_{k=0}^{d}\EE \Bigg[\sum_{ \underset { (y_1, \ldots, y_k)_{\neq}\in \mathbf{X}^k}{ (x_1, \ldots, x_{d+1})_{\neq} \in \mathbf{X}^{d+1}}} \mathbb{1}_{f(\mathbf{x}_{d+1}) > v_{\rho}}\mathbb{1}_{f(\mathbf{x}_{1:k}, \mathbf{y}_{k+1:d+1})>v_{\rho}}\mathbb{1}_{R(\mathbf{x}_{1:d+1})\geq R(\mathbf{x}_{1:k}, \mathbf{y}_{k+1:d+1})}\\
\times \mathbb{1}_{\mathbf{X}\cap B^{(\cup)}(\mathbf{x}_{1:d+1}, \mathbf{y}_{k+1:d+1}) - \{\mathbf{x}_{1:d+1}\}\cup \{ \mathbf{y}_{k+1:d+1}\}} = \varnothing \Bigg].     
\end{multline*}
where $B^{(\cup)}(\mathbf{x}_{1:d+1}, \mathbf{y}_{k+1:d+1})=B(\mathbf{x}_{1:d+1})\cup B(\mathbf{x}_{1:k}, \mathbf{y}_{k+1:d+1} )$. It results of Slivnyak's formula (see e.g. Theorem 3.3.5 in \cite{SW}) that 
\begin{multline*}
\PPP{f(\cell)>v_{\rho}} = 2\sum_{k=0}^{d}\int_{(\RR^d)^{d+1-k}}\int_{(\RR^d)^{d+1}}\mathbb{1}_{f(\mathbf{x}_{d+1}) > v_{\rho}}\mathbb{1}_{f(\mathbf{x}_{1:k}, \mathbf{y}_{k+1:d+1})>v_{\rho}}\mathbb{1}_{R(\mathbf{x}_{1:d+1})\geq R(\mathbf{x}_{1:k}, \mathbf{y}_{k+1:d+1})}\\
\times \PPP{\#\mathbf{X}\cap B^{(\cup)}(\mathbf{x}_{1:d+1}, \mathbf{y}_{k+1:d+1}) = 0 }d\mathbf{x}_{1:d+1} d\mathbf{y}_{k+1:d+1}.
\end{multline*}
We conclude the proof of Lemma \ref{genericdelaunay} noting that $\#\mathbf{X}\cap B^{(\cup)}(\mathbf{x}_{1:d+1}, \mathbf{y}_{k+1:d+1})$ is Poisson distributed of mean $\beta_d\lambda_d\left(B^{(\cup)}(\mathbf{x}_{1:d+1}, \mathbf{y}_{k+1:d+1}) \right)$ and using for all  $y_{k+1}, \ldots, y_{d+1}$ the (Blaschke-Petkantschin type) change of variables \begin{equation}\label{Blaschke} \begin{split} \phi_1: & \RR_+\times \RR^d \times (\mathbf{S}^{d-1})^{d+1}\longrightarrow (\RR^d)^{d+1}\\
& (r,z,\mathbf{u}_{1:d+1}) \longmapsto \mathbf{x}_{1:d+1} \text{ with } x_i = z+ru_i\end{split}\end{equation} where the Jacobian matrix is given by $|D\phi_1(r,z,\mathbf{u}_{1:d+1})| = r^{d^2-1}\lambda_d(\mathbf{u}_{1:d+1})$.

\end{prooft}

In Lemma \ref{genericdelaunay}, we have assumed that $R\left(r\mathbf{u}_{1:k}, \mathbf{y}_{k+1:d+1}\right)$ is less than   $R(r\mathbf{u}_{1:d+1})$. It overcomes the difficulty to consider elongated cells. This property will be needed in sections \ref{maxareasection} and \ref{minareasection} but not in section \ref{mincirconssection} since we consider small circumradii.

\subsection{Minimum of the circumradii}
\label{mincirconssection}
Let us recall that $R(C)$ denotes the circumradius of the cell $C\in\mm_{PDT}$. In this subsection, we investigate the minimum 
\[R_{\min, PDT}(\rho) = \min_{\underset{z(C)\in \mathbf{W}_{\rho}}{C\in\mm_{PDT}},}R(C).\] The asymptotic behaviour of $R_{\min, PDT}(\rho)$ is given in the following proposition.

\begin{Prop}
\label{mincirconsdelappl}
Let $\mm_{PDT}$ be a Poisson-Delaunay tessellation of intensity $\gamma = 1$ in $\RR^d$, $d\geq 2$. Then for all $t\geq 0$
\begin{equation}\label{mincirconsdel}\left|\PPP{\alpha_{d,1}^{1/d}\rho^{1/d}R_{\min, PDT}(\rho)^d\geq t} - e^{-t^{d}}\right| = O\left(\rho^{-1/d} \right)\end{equation} where 
\[\alpha_{d,1} = \frac{\delta_d^d}{d!} = \frac{(\kappa_d\beta_d)^d}{d!} = \frac{1}{d!}\cdot \left(\frac{(d^3+d^2)\Gamma\left(\frac{d^2}{2} \right) \Gamma^d\left(\frac{d+1}{2} \right)\pi^{1/2}}{2^{d+1}\Gamma\left(\frac{d^2+1}{2} \right)\Gamma^{d+1}\left(\frac{d+2}{2} \right)}\right)^d  .\]
\end{Prop}
The asymptotic behaviour of the maximum of circumradii has been investigated in \cite{CC} and will be recalled in section \ref{sectionrealtyp2}.

\begin{prooft}{Proposition \ref{mincirconsdelappl}}
First, we give the asymptotic behaviour of the distribution function of $R(\cell)$. According to  \eqref{typicalcelldelaunay}, the random variable $R(\cell)^d$ is Gamma distributed of parameters $\left(d^2,\delta_d^{-1}\right)$. Thanks to consecutive integration by parts, this provides that  \begin{equation}\label{circonsdel}\PP(R(\cell)<v) = \sum_{i=d}^{\infty}\frac{1}{i!}{(\delta_dv^{d})}^ie^{-\delta_dv^d}\end{equation} for all $v\geq 0$. A Taylor approximation of the right-hand side  when $v$ is small shows that $|\PP(R(\cell)<v) - \alpha_{d,1}\cdot v^{d^2}|$ is of order $v^{d^2+d}$. Hence, taking for all $t\geq 0$
\begin{equation}
\label{defvmincirconsdel}
v_{\rho} = v_{\rho}(t) = \left( \alpha_{d,1}^{-1}\rho^{-1}\right)^{1/d^2}t^{1/d}
\end{equation} 
we obtain
\begin{equation}
\label{G1mincirconsdel}
G_1(\rho) = |\rho \PP(R(\cell)<v_{\rho}) - t^d  | = O\left(\rho^{-1/d} \right).
\end{equation}

To calculate the order of $G_2(\rho)$, it is enough to give a suitable upper bound of $G_{2,k}(\rho)$ for all $k=0,\ldots, d$ according to Lemma \ref{genericdelaunay}. Bounding the exponential in \eqref{defg2kdel} by 1 (a suitable estimate when considering small cells) and $\lambda_d(\mathbf{u}_{1:d+1})$ by a constant, we deduce for all $r\in\RR_+$, $\mathbf{u}_{1:d+1}\in (\mathbf{S}^{d-1})^{d+1}$ and $\mathbf{y}_{k+1:d+1}\in(\RR^d)^{d+1-k}$ that 
\begin{equation}
\label{g2kmincirconsdel}
g_{2,k}(\rho,r,\mathbf{u}_{1:d+1},\mathbf{y}_{k+1:d+1}) \leq c\cdot r^{d^2-1}\mathbb{1}_{r<v_{\rho}} \mathbb{1}_{E_{k,r,\mathbf{u}_{1:d+1}}}(\mathbf{y}_{k+1:d+1}).
\end{equation}  

When $k=0$, we bound $\mathbb{1}_{E_{0,r,\mathbf{u}_{1:d+1}}}(\mathbf{y}_{1:d+1})$ by $\mathbb{1}_{R( \mathbf{y}_{1:d+1})<v_{\rho}}\cdot\mathbb{1}_{z(\mathbf{y}_{1:d+1})\in\CC}$. We can omit the last condition in \eqref{delcond1} and the two conditions in \eqref{delcond2} since having a small circumradius almost guarantees that they are satisfied. Integrating the right-hand side of \eqref{g2kmincirconsdel} and taking the same change of variables as in \eqref{Blaschke} i.e. $y_i=z'+r'u'_i$, $i=1,\ldots, d+1$, we deduce from \eqref{defG2k} and \eqref{defvmincirconsdel} that 
 \begin{equation}
\label{g20mincirconsdel}
G_{2,0}(\rho)\leq c\cdot\rho\int_0^{v_{\rho}}r^{d^2-1}dr\times \lambda_d(\CC)\int_{0}^{v_{\rho}}r'^{d^2-1}dr' = O\left(\log\rho\cdot \rho^{-1} \right).
\end{equation}
 
When $k=1,\ldots, d$, we use the fact that $R(r\mathbf{u}_{1:k},\mathbf{y}_{k+1:d+1})<v_{\rho}\Longrightarrow y_i\in B(ru_1,2v_{\rho})$ for all $i=k+1,\ldots, d+1$. Bounding $\mathbb{1}_{E_{k,r,\mathbf{u}_{1:d+1}}}(\mathbf{y}_{k+1:d+1})$ by $\mathbb{1}_{y_{k+1}, \ldots, y_{d+1}\in B(ru_1,2v_{\rho})}$ and integrating \eqref{g2kmincirconsdel}, we deduce from \eqref{defG2k} that
\begin{multline} 
\label{majG2kmincirconsdel}
G_{2,k}(\rho) \leq c\cdot \rho\int_0^{v_{\rho}} \int_{\mathbf{S}^{d-1}}  \int_{(\RR^d)^{d+1}}r^{d^2-1}\mathbb{1}_{y_{k+1}, \ldots, y_{d+1}\in B(ru_1,2v_{\rho})}drd\sigma(u_1)d\mathbf{y}_{k+1:d+1}\\  \leq c\cdot \rho\int_0^{v_{\rho}}r^{d^2-1}dr\times v_{\rho}^{d(d+1-k)} = O\left(\rho^{-(d+1-k)/d} \right).\end{multline} 
Since $k=0,\ldots, d$, the right-hand side of \eqref{majG2kmincirconsdel} is less than $\rho^{-1/d}$ for $\rho$ large enough. Indeed, $G_{2,k}(\rho)$ is maximal when $k=d$ i.e. when the two distinct Delaunay cells have $d$ common vertices.  From \eqref{majG2del}, \eqref{g20mincirconsdel} and \eqref{majG2kmincirconsdel} we deduce that \begin{equation}\label{mincirconsdeltyp2}G_2(\rho) = O\left(\rho^{-1/d} \right).\end{equation}
The rate of convergence \eqref{mincirconsdel} is now a direct consequence of \eqref{G1mincirconsdel},  \eqref{mincirconsdeltyp2} and Theorem \ref{realtyp}. 
\end{prooft}

When $d=1$, the order of $R_{\min, PDT}(\rho)$ is $\rho^{-1}$. Moreover, the rate of convergence is $\log\rho\cdot \rho^{-1}$ (and not $\rho^{-1}$) since this is the order of  $\PP(A_{\rho})$ and $N_{\rho}^{-1}$ which appear in Theorem \ref{realtyp}. 

Let us remark that a slightly weaker version of Proposition \ref{mincirconsdelappl} in $\RR^d$ could have been deduced from a theorem due Schulte and Thäle (see Theorem 1.1 in \cite{ST}). It comes from the fact that $R_{\min, PDT}(\rho)$ can be written as a minimum of a $U$-statistic. More precisely
\[R_{\min, PDT}(\rho) = \min_{\underset{z\left(\mathbf{x}_{1:d+1}\right)\in \mathbf{W}_{\rho}}{\mathbf{x}_{1:d+1}\in\mathbf{X}^{d+1}},}R(\mathbf{x}_{1:d+1}).\] Indeed, if a simplex induced by a set of $(d+1)$ distinct points $\mathbf{x}_{1:d+1}$ of $\mathbf{X}$ minimizes the circumradius, it is necessarily a Delaunay cell: otherwise, the circumball $B(\mathbf{x}_{1:d+1})$ contains a point of $\mathbf{X}$ in its interior which contradicts the minimality of $R(\mathbf{x}_{1:d+1})$. Nevertheless, the rate of convergence $O\left(\rho^{-1/d}\right)$ of Proposition \ref{mincirconsdelappl} is more accurate than the rate deduced from Theorem 1.1. in \cite{ST} since the latter is of order $O\left(\rho^{-1/2d}\right)$. To the best of our knowledge, the convergence of the point process provided by Theorem \ref{PPTh} applied to the circumscribed radius of Delaunay cells is new.

\subsection{Maximum of the areas, $d=2$}
\label{maxareasection}
Here and in the subsequent subsection, we investigate the extremes of the areas of a planar Poisson-Delaunay tessellation of intensity 1. The extension to higher dimension would be intricate since the integral formula for the distribution function of the volume of the typical cell becomes intractable. The intensity of the underlying Poisson point process is \begin{equation}\label{defbeta}\gamma_{\mathbf{X}} = \beta_2 = \frac{1}{2}.\end{equation}
In this subsection, we investigate the maximum of the areas i.e. \[A_{\max, PDT}(\rho)=\max_{\underset{z(C)\in \mathbf{W}_{\rho}}{C\in\mm_{PDT}},}\lambda_2(C).\] 
The following proposition shows that $A_{\max, PDT}(\rho)$ is of order $\log\rho$.

\begin{Prop}
\label{maxvoldelappl}
Let $\mm_{PDT}$ be a Poisson-Delaunay tessellation of intensity $\gamma=1$ in $\RR^2$. Then for all $t\in\RR$
\begin{equation}\label{maxareadeltess}\left|\PPP{\alpha_{2} A_{\max, PDT}(\rho)-\log\left(\frac{3}{2}\rho \right)\leq t} - e^{-e^{-t}}\right| = O\left(1/\log\rho\right)\end{equation} where
\begin{equation}\label{defalpha45}\alpha_{2}=\frac{2\pi}{3\sqrt{3}}.\end{equation}
\end{Prop}

\begin{prooft}{Proposition \ref{maxvoldelappl}}
Thanks to \eqref{typicalcelldelaunay}, the distribution function of $\lambda_2(\cell)$ can be made explicit. Indeed, an integral representation of $\PPP{\lambda_2(\cell)>v}$ due to Rathie (see (3.2) in \cite{Ra}) is 
\begin{equation}\label{densityvol}\PP(\lambda_2(\cell)>v)=\frac{6}{\pi}\int_{\alpha_{2}\beta_2v}^{\infty}xK^2_{1/6}(x)dx \end{equation} 
where $K_{1/6}(\cdot)$ denotes the modified Bessel function of order $1/6$. When $x$ goes to infinity, a Taylor approximation of $K_{1/6}(x)$ is given by (see Formula 9.7.2, p.  378 in \cite{AS})
\begin{equation}
\label{DLmaxvoldel}
K_{1/6}(x) = \sqrt{\frac{\pi}{2x}}e^{-x}\left(1+O\left(\frac{1}{x}\right) \right).
\end{equation}
We deduce from \eqref{defbeta}, \eqref{densityvol} and \eqref{DLmaxvoldel} that for $v$ large enough
\begin{equation}
\label{expt}
\left|\PPP{\lambda_2(\cell)>v} - \frac{3}{2} e^{-\alpha_{2}v}\right|
 \leq c\cdot \int_{\frac{1}{2}\alpha_{2}v}^{\infty}\frac{e^{-2x}}{x}dx\leq c\cdot \frac{e^{-\alpha_{2}v}}{v}.
\end{equation}
Taking for all $t\in\RR$
\begin{equation}
\label{defvmaxvol}
v_{\rho} = v_{\rho}(t) = \frac{1}{\alpha_{2}}\left(\log\left(\frac{3}{2}\rho \right) +t\right).
\end{equation}
we obtain from \eqref{expt} that
\begin{equation}
\label{maxvoldeltyp}
G_1(\rho) = |\rho\PPP{\lambda_2(\cell)>v_{\rho}} - e^{-t}| = O\left(1/\log\rho\right).
\end{equation}

In the rest of the proof, we give a suitable upper bound of $G_2(\rho)$. Taking $f(\cdot) = \lambda_2(\cdot)$ in \eqref{defg2kdel} and using the facts that $\lambda_2(r\mathbf{u}_{1:3}) = r^2\lambda_2(\mathbf{u}_{1:3})$ and   $\lambda_2(\mathbf{u}_{1:3})\leq c$, we have

\begin{equation}
\label{majG2delmaxvol}
g_{2,k}(\rho,r,\mathbf{u}_{1:3},\mathbf{y}_{k+1:3}) \leq c\cdot r^{3}e^{-\frac{1}{2}\lambda^{(\cup)}_d(r, \mathbf{u}_{1:k}, \mathbf{y}_{k+1:3})} \mathbb{1}_{r^2\lambda_2(\mathbf{u}_{1:3})>v_{\rho}}\mathbb{1}_{E_{k,r,\mathbf{u}_{1:3}}}(\mathbf{y}_{k+1:3})  .
\end{equation} for all $k=0,1,2$. To bound $g_{2,k}(\cdot)$, the key idea is to give a suitable lower bound of the area of the union of two disks (see Figure \ref{figarea} (a)). This is provided in the following fundamental lemma.

\begin{Le}
\label{trianglelemma}
Let $\{\mathbf{x}_{1:3}\}=\{x_1,x_2,x_3\}$ and $\{\mathbf{x}'_{1:3}\}=\{x'_1,x'_2,x'_3\}$ be two 3-tuples of points in $\RR^2$ such that $x_i\not\in B(\mathbf{x}'_{1:3})$ and $x'_j\not\in B(\mathbf{x}_{1:3})$ for all $i,j=1,2,3$. Let us assume that $R:=R(\mathbf{x}_{1:3})\geq R(\mathbf{x}'_{1:3})$. Then
 \begin{equation}\label{triangle}
\lambda_2\left(B(\mathbf{x}_{1:3})\cup B(\mathbf{x}'_{1:3})\right) \geq  \left(\frac{\pi}{2}-1\right)R^2 + \lambda_2(\mathbf{x}_{1:3}) + \lambda_2(\mathbf{x}'_{1:3}).
\end{equation} \end{Le}

\begin{prooft}{Lemma \ref{trianglelemma}}
Let $\{\mathbf{x}_{1:3}\}$ and $\{\mathbf{x}'_{1:3}\}$ be two 3-tuples in $\RR^2$.

If the interior of $B(\mathbf{x}_{1:3})\cap B(\mathbf{x}'_{1:3})$ is empty, we have
\begin{equation}\label{minvolunion}\lambda_2(B(\mathbf{x}_{1:3})\cup B(\mathbf{x}'_{1:3})) = \lambda_2(B(\mathbf{x}_{1:3})) + \lambda_2(B(\mathbf{x}'_{1:3}))\geq \pi R^2+\lambda_2(\mathbf{x}'_{1:3}).\end{equation} Moreover, the maximal area of a triangle inscribed in a ball of radius $R$ is $\frac{3\sqrt{3}}{4}R^2$ which is the area of an equilateral triangle. In particular, we have $\lambda_2(\mathbf{x}_{1:3})\leq \frac{3\sqrt{3}}{4}R^2$. This together with \eqref{minvolunion} implies that
\[\lambda_2(B(\mathbf{x}_{1:3})\cup B(\mathbf{x}'_{1:3})) \geq \left(\pi-\frac{3\sqrt{3}}{4} \right)R^2+\lambda_2(\mathbf{x}_{1:3})+ \lambda_2(\mathbf{x}'_{1:3})\geq \left( \frac{\pi}{2}-1\right)R^2+\lambda_2(\mathbf{x}_{1:3})+ \lambda_2(\mathbf{x}'_{1:3}) \]

If $B(\mathbf{x}_{1:3})\cap B(\mathbf{x}'_{1:3})$ has non empty interior, the intersection of the circumspheres induced by the points $\mathbf{x}_{1:3}$ and $\mathbf{x}'_{1:3}$ is reduced to two points, say $p_1, p_2\in\RR^2$. Let us denote by $\mathbf{L}$ the affine line $(p_1,p_2)$ and $\mathbf{H}^-$ (respectively $\mathbf{H}^+$) the half plane delimited by $\mathbf{L}$ and containing (respectively not containing) the circumcenter $z(\mathbf{x}_{1:3})$. Since $x_i\not\in B(\mathbf{x}'_{1:3})$ and $x'_j\not\in B(\mathbf{x}_{1:3})$, $i,j=1,2,3$, the triangle $\Delta(\mathbf{x}'_{1:3})$ is included in $\mathbf{H}^+$. Hence
 \begin{multline}\label{ballincl}\lambda_2\left(B(\mathbf{x}_{1:3})\cup B(\mathbf{x}'_{1:3})\right) = \lambda_2\left((B(\mathbf{x}_{1:3})\cup B(\mathbf{x}'_{1:3}))\cap \mathbf{H}^-\right) + \lambda_2\left((B(\mathbf{x}_{1:3})\cup B(\mathbf{x}'_{1:3}))\cap \mathbf{H}^+\right)\\
  \geq  \lambda_2(B(\mathbf{x}_{1:3})\cap \mathbf{H}^-) + \lambda_2(\mathbf{x}'_{1:3}).\end{multline}

\begin{figure}\begin{center}\begin{tabular}{cc}\includegraphics[width=8cm,height=6.5cm]{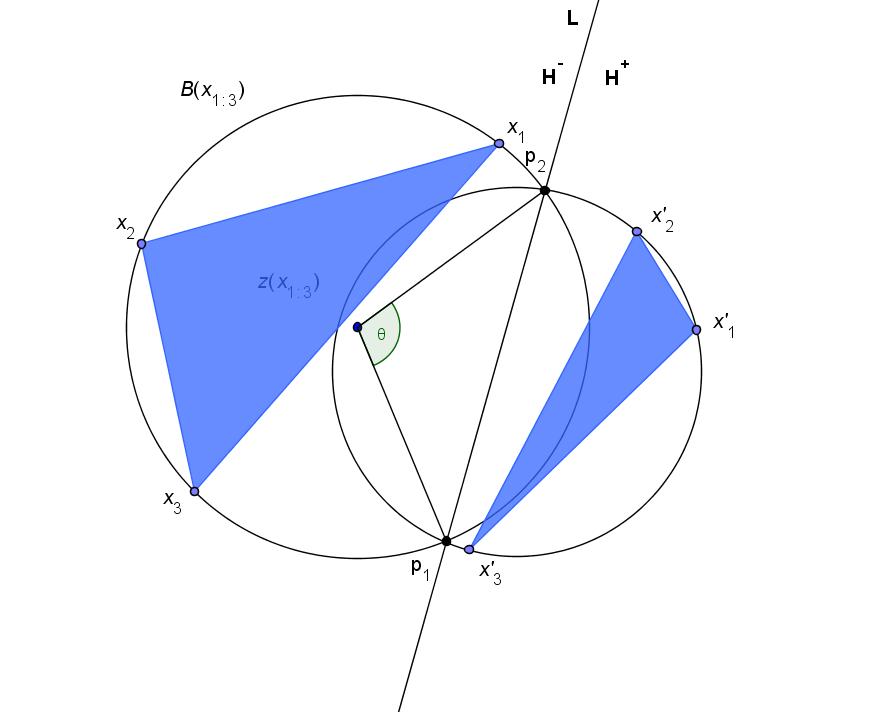} & \includegraphics[width=8cm,height=6.5cm]{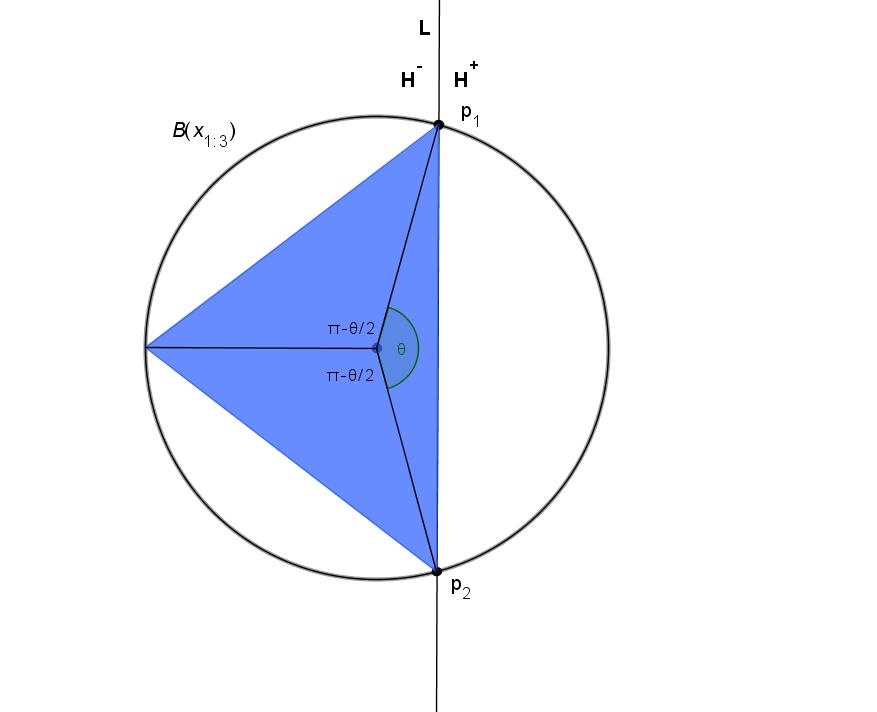}  \\ (a) & (b)\end{tabular}\end{center}\caption{(a). A union of two disks. (b). The triangle which maximizes the area.} \label{figarea}\end{figure}


 In the rest of the  proof, we provide a suitable lower bound of $ \lambda_2(B(\mathbf{x}_{1:3})\cap \mathbf{H}^-)$. To do it, we denote by $\theta\in [0,2\pi]$ the angle  $\angle p_1z(\mathbf{x}_{1:3})p_2$. Actually $\theta\in [0,\pi]$: this comes from the fact that  $\lambda_2(B(\mathbf{x}_{1:3})\cap \mathbf{H}^-)\geq\frac{\pi}{2}R^2$ since $R:=R(\mathbf{x}_{1:3})\geq R(\mathbf{x}'_{1:3})$. The area of the cap $B(\mathbf{x}_{1:3})\cap \mathbf{H}^-$ is given by
 \begin{equation}
 \label{exprcap}
 \lambda_2(B(\mathbf{x}_{1:3})\cap \mathbf{H}^-) = \left(\pi - \frac{1}{2}(\theta-\sin\theta)\right)R^2.
 \end{equation} We discuss below two cases depending on $\theta$.
 
 If $\theta\in [0,2\pi/3]$, we deduce from \eqref{exprcap} that \begin{equation}\label{mincap1}
  \lambda_2(B(\mathbf{x}_{1:3})\cap \mathbf{H}^-) \geq \left(\frac{2\pi}{3}+\frac{\sqrt{3}}{4} \right)R^2.
 \end{equation} Since $\lambda_2(\mathbf{x}_{1:3})$ is less than $\frac{3\sqrt{3}}{4}R^2$, we deduce from \eqref{mincap1} that \begin{equation}
 \label{mincap2}
  \lambda_2(B(\mathbf{x}_{1:3})\cap \mathbf{H}^-) \geq \lambda_2(\mathbf{x}_{1:3}) + \left(\frac{2\pi}{3}-\frac{\sqrt{3}}{2} \right)R^2\geq \lambda_2(\mathbf{x}_{1:3}) + \left(\frac{\pi}{2}-1 \right)R^2.
 \end{equation} In that case, the inequality \eqref{triangle} results from \eqref{ballincl} and \eqref{mincap2}.
 
 If $\theta\in [2\pi/3, \pi]$, with a standard method of geometry, we can show that the maximal area of a triangle inscribed in $B(\mathbf{x}_{1:3})\cap \mathbf{H}^-$, denoted by $M(\theta)$, is \begin{equation}\label{defMaxtriangle} M(\theta) = \left(\sin\frac{\theta}{2} + \frac{1}{2}\sin\theta \right)R^2.\end{equation} Actually, the triangle which maximizes the area is isoscele with central angles $\pi-\theta/2, \pi-\theta/2$ and $\theta$ (see Figure \ref{figarea} (b)). In particular, we have \begin{equation}\label{regularsimplexmaj}
 \lambda_2(\mathbf{x}_{1:3})\leq M(\theta).\end{equation}  We obtain from \eqref{exprcap} and \eqref{defMaxtriangle} that \begin{equation}
  \lambda_2(B(\mathbf{x}_{1:3})\cap \mathbf{H}^-) \geq M(\theta) + \left(\frac{\pi}{2}-1 \right)R^2 + \left(\frac{\pi}{2}+1-\left(\frac{1}{2}\theta + \sin\frac{\theta}{2} \right) \right)R^2.
 \end{equation} The last term of the right-hand side is a decreasing function on $[0,\pi]$. Its minimum equals 0 at $\theta=\pi$ i.e.
 \begin{equation*}
 \frac{\pi}{2}+1-\left(\frac{1}{2}\theta + \sin\frac{\theta}{2} \right) \geq 0
 \end{equation*} for all $\theta\in [0,\pi]$. This shows that 
 \begin{equation}\label{mincap4}\lambda_2(B(\mathbf{x}_{1:3})\cap \mathbf{H}^-) \geq M(\theta) + \left(\frac{\pi}{2}-1 \right)R^2.\end{equation}
  The inequality \eqref{triangle} is a direct consequence of  \eqref{ballincl}, \eqref{regularsimplexmaj} and \eqref{mincap4}. 
 
   \end{prooft}

We can now derive an upper bound of $g_{2,k}(\cdot)$ for all $k=0,1,2$. Indeed, if $\mathbf{y}_{k+1:3}\in E_{k,r,\mathbf{u}_{1:3}}$, where $E_{k,r,\mathbf{u}_{1:3}}$ has been defined in \eqref{defE}, the set of points $\{\mathbf{x}_{1:3}\} = \{r\mathbf{u}_{1:3}\}$ and $\{\mathbf{x}'_{1:3}\}=\{r\mathbf{u}_{1:k},\mathbf{y}_{k+1:3}\}$ satisfies the assumptions of Lemma \ref{trianglelemma} since  $R(r\mathbf{u}_{1:3}) = r$ and $R(r\mathbf{u}_{1:3})\geq R(r\mathbf{u}_{1:k},\mathbf{y}_{k+1:3})$. Using the fact that $B(r\mathbf{u}_{1:3}) = B(0,r)$, $\lambda_2(r\mathbf{u}_{1:3})>v_{\rho}$ and $\lambda_2(r\mathbf{u}_{1:k}, \mathbf{y}_{k+1:3})>v_{\rho}$, we deduce from \eqref{defunionvol}, \eqref{majG2delmaxvol} and \eqref{triangle} that

\begin{equation}
\label{majG2delmaxvol2} 
g_{2,k}(\rho, r, \mathbf{u}_{1:3},\mathbf{y}_{k+1:3}) \leq c\cdot r^{3}e^{-\frac{1}{2}\left( \left(\frac{\pi}{2}-1\right) r^2 + 2v_{\rho}  \right)}
 \mathbb{1}_{r^2\lambda_2(\mathbf{u}_{1:3})>v_{\rho}} \mathbb{1}_{E_{k,r,\mathbf{u}_{1:3}}}(\mathbf{y}_{k+1:3}).
\end{equation} 

Since $\frac{3\sqrt{3}}{4}r^2\geq r^2\lambda_2(\mathbf{u}_{1:3})$, we deduce from \eqref{defalpha45} and \eqref{defvmaxvol}  that \begin{equation}\label{implicationindicator} r^2\lambda_2(\mathbf{u}_{1:3})>v_{\rho} \Longrightarrow r^2>4v_{\rho}/3\sqrt{3} \Longrightarrow r> (2\left(\log\rho + c \right)/\pi)^{1/2} \end{equation} where $c=\log(3/2)+t$. Integrating the right-hand side on $\mathbf{y}_{k+1:3}$, we obtain

\begin{equation}
\label{majG2delmaxvol3}
G_{2,k}(\rho) \leq c\cdot\rho\int_{(2\left(\log\rho + c \right)/\pi)^{1/2}}^{\infty}\int_{(\mathbf{S}^1)^3}r^{3}e^{-\frac{1}{2}\left( \left(\frac{\pi}{2}-1\right) r^2 + 2v_{\rho}  \right)}\times \lambda_{2(3-k)}(E_{k,r, \mathbf{u}_{1:3}})drd\sigma(\mathbf{u}_{1:3}).
\end{equation} The following lemma gives a uniform upper bound of  $\lambda_{2(3-k)}(E_{k,r, \mathbf{u}_{1:3}})$. 

\begin{Le}
\label{lebesgueElemma}
Let $\mathbf{u}_{1:3}\in (\mathbf{S}^1)^k$ and $r>(2\left(\log\rho + c \right)/\pi)^{1/2}$. Then for $\rho$ large enough
\begin{equation}
\label{lebesgueE}
\lambda_{2(3-k)}(E_{k,r, \mathbf{u}_{1:3}}) \leq c\cdot r^{2(3-k)}.
\end{equation}
\end{Le}

\begin{prooft}{Lemma \ref{lebesgueElemma}}
We discuss three cases that depend on $k$. 

If $k=2$, we show that $E_{2,r, \mathbf{u}_{1:3}}$ is included in a ball of radius $r$ up to a multiplicative constant and centered at 0. Let $y_3$ be in $E_{2,r, \mathbf{u}_{1:3}}$. From the triangle inequality, we have
\begin{equation}\label{lebesgue1}|y_3| \leq |y_3 - z(r\mathbf{u}_{1:2},y_3)| + |z(r\mathbf{u}_{1:2}, y_3)| \leq  r+\text{diam}(\CC).\end{equation} The last inequality comes from the fact that $|y_3 - z(r\mathbf{u}_{1:2},y_3)|$ is the circumradius of $\Delta(r\mathbf{u}_{1:2}, y_3)$, which is less than $r$, and the fact that $z(r\mathbf{u}_{1:2})\in\CC$. Moreover
\begin{equation}\label{lebesgue2}\text{diam}(\CC) \leq c\cdot (\log\rho)^{1/2}\leq c\cdot r\end{equation} where the last inequality holds for $\rho$ large enough since $r>(2\left(\log\rho + c \right)/\pi)^{1/2}$ converges to $\infty$ as $\rho$ goes to infinity. We deduce from \eqref{lebesgue1} and \eqref{lebesgue2} that 
\begin{equation}\label{lebesgue3}|y_3| \leq c\cdot r\end{equation}  The upper bound \eqref{lebesgue3} shows that $E_{2,r, \mathbf{u}_{1:3}}\subset B(0,c\cdot r)$. In particular, 
\[\lambda_2(E_{2,r, \mathbf{u}_{1:3}})\leq c\cdot r^2.\] 

If $k=1$ or $k=0$, proceeding along the same lines as in the case $k=2$, we show that $E_{k,r, \mathbf{u}_{1:3}}\subset B(0,c\cdot r)^{3-k}$ and consequently we get $\lambda_{2(3-k)}(E_{k,r, \mathbf{u}_{1:3}})\leq c\cdot r^{2(3-k)}.$

\end{prooft}

We can now derive an upper bound of $G_{2,k}(\rho)$. Indeed, integrating $\mathbf{u}_{1:3}$ on $(\mathbf{S}^1)^3$, we deduce from \eqref{majG2delmaxvol3} and \eqref{lebesgueE} that
\[ G_{2,k}(\rho)\leq c\cdot\rho\int_{(2\left(\log\rho + c \right)/\pi)^{1/2}}^{\infty}r^{9-k}e^{-\frac{1}{2}\left( \left(\frac{\pi}{2}-1\right) r^2 + 2v_{\rho}  \right)}dr.\] Integrating the right-hand side, we obtain from \eqref{defvmaxvol} that 
\begin{equation} \label{ratemaxvol} G_{2,k}(\rho) \leq c\cdot (\log\rho)^{8-2k}\rho^{(\pi+2-3\sqrt{3})/2\pi}=O\left((\log\rho)^8\rho^{-\epsilon}\right)\end{equation} with $\epsilon = -\pi-2+3\sqrt{3}>0$. Proposition 2 results of \eqref{ratemaxvol}, Lemma \ref{genericdelaunay} and Theorem \ref{realtyp}.
\end{prooft}

Lemma \ref{trianglelemma} provides the main tool of the proof. We can note that the inequality \eqref{triangle} is obvious when we replace $\frac{\pi}{2}-1$ by a constant $\alpha\leq \pi-\frac{3\sqrt{3}}{2}$. Indeed, if $\Delta(\mathbf{x}_{1:3})$ and $\Delta(\mathbf{x}'_{1:3})$ are two triangles with $R:=R(\mathbf{x}_{1:3})\geq R(\mathbf{x}'_{1:3})$, a trivial inequality is
\[\lambda_2(B(\mathbf{x}_{1:3})\cup B(\mathbf{x}'_{1:3}))\geq \pi R^2.\]
Consequently \[\lambda_2(B(\mathbf{x}_{1:3})\cup B(\mathbf{x}'_{1:3}))\geq \left(\pi-\frac{3\sqrt{3}}{2} \right)R^2+\lambda_2(\mathbf{x}_{1:3})+\lambda_2(\mathbf{x}'_{1:3})\] since $\lambda_2(\mathbf{x}_{1:3})$ and $\lambda_2(\mathbf{x}'_{1:3})$ are less than $\frac{3\sqrt{3}}{4}R^2$. Nevertheless, the previous lower bound is not enough to guarantee that $G_{2,k}(\rho)$ converges to 0. The important fact in Lemma \ref{trianglelemma} is that we consider the more precise constant $\frac{\pi}{2}-1>\pi-\frac{3\sqrt{3}}{2}$.

Another remark deals with the shape of the cell maximizing the area. As we will see in Example 2 of section \ref{sectionrealtyp2}, the maximum of circumradii of a planar Poisson-Delaunay tessellation, denoted by $R_{\max, PDT}(\rho)$, is of order $(\delta_2^{-1}\log\rho)^{1/2} = (2\pi^{-1}\log\rho)^{1/2}$ according to \eqref{defdelta} and \eqref{maxcircumdel}. Thanks to \eqref{maxareadeltess}, this shows that $A_{\max, PDT}(\rho)$ equals asymptotically $\frac{3\sqrt{3}}{4}R_{\max, PDT}^2(\rho)$ which is the area of an equilateral triangle of circumradius $R_{\max, PDT}(\rho)$. It seems that the shape of the cell maximizing the area tends to that of an equilateral triangle. This fact can be connected to the D.G. Kendall's conjecture and to the work of Hug and Schneider in \cite{HS}.

\subsection{Minimum of the areas, $d=2$}
\label{minareasection}
In our third example, we calculate the asymptotic behaviour of the minimum of the areas of the cells of a Poisson-Delaunay tessellation (of intensity 1) in $\RR^2$ i.e. \[A_{\min, PDT}(\rho)=\min_{\underset{z(C)\in \mathbf{W}_{\rho}}{C\in\mm_{PDT}},}\lambda_2(C).\] The asymptotic behaviour is given in the following proposition.

\begin{Prop}
\label{minvoldelappl}
Let $\mm_{PDT}$ be a Poisson-Delaunay tessellation of intensity $\gamma=1$ in $\RR^2$. Then for all $t\geq 0$
\begin{equation}\label{minvoldel} \PPP{\alpha_{3}^{3/5}\rho^{3/5}A_{\min,PDT}(\rho)\geq t} \conv[\rho]{\infty}e^{-t^{5/3}} \end{equation} where \[\alpha_{3} = 2^{-2/3}\cdot 3^{-1/2}\cdot 5^{-1}\cdot \pi^{2/3}\cdot \Gamma(1/6)^2.\] 
\end{Prop}

In \cite{ST}, Schulte and Thäle investigate the behaviour of the smallest area $S_{\rho}$  of all triangles that can be formed by three points of the Poisson point process i.e. 
\[S_{\rho} = \min_{\underset{z\left(\mathbf{x}_{1:3}\right)\in \mathbf{W}_{\rho}}{\mathbf{x}_{1:3}\in\mathbf{X}^{3}},}\lambda_2(\mathbf{x}_{1:3}).\] The asymptotic behaviour of $S_{\rho}$ is given by (see Theorem 2.5. in \cite{ST})
\[\PPP{\rho S_{\rho}\geq t} \conv[\rho]{\infty}e^{-\beta t}\] where $\beta$ is a constant which can be made explicit. The previous limit compared to \eqref{minvoldel} shows that the smallest area of the Delaunay cells is much larger than the smallest area of all triangles. 

\begin{prooft}{Proposition \ref{minvoldelappl}}
First, we calculate the asymptotic behaviour of the distribution function of $\lambda_2(\cell)$. Let us recall that such a function is given in \eqref{densityvol}. A Taylor expansion of the modified Bessel function of order 1/6 is given by (see Formula 9.6.9, p.  375 in \cite{AS}) \begin{equation}\label{DLminvoldel}K_{1/6}(x) = 2^{-5/6}\Gamma\left(1/6 \right)x^{-1/6} + o(x^{-1/6}).\end{equation} This together with \eqref{defbeta}  and \eqref{densityvol} shows that
\begin{equation}\label{Taylorminvol}\PPP{\lambda_2(\cell)<v} = \frac{6}{\pi}\cdot 2^{-5/3}\Gamma\left(1/6\right)^2\int_0^{\alpha_{2}\beta_2 v}\left(x^{2/3}+o\left(x^{2/3}\right)\right)dx = \alpha_{3}\cdot v^{5/3} + o\left(v^{5/3}\right).\end{equation}  Taking for all $t\geq 0$
\begin{equation}
\label{defvminvol}
v_{\rho} = v_{\rho}(t) = (\alpha_{3}^{-1}\rho^{-1})^{3/5}t
\end{equation}
we obtain
\begin{equation}
\label{majG1minvol}
G_1(\rho) = |\rho\PPP{\lambda_2(\cell)<v_{\rho}} - t^{5/3}| \conv[\rho]{\infty}0.
\end{equation}

We investigate below the rate of convergence of $G_2(\rho)$. Taking $f(r\mathbf{u}_{1:3})= r^2\lambda_2(\mathbf{u}_{1:3})$ and using the fact that $\lambda_2(B(0,r)\cup B(r\mathbf{u}_{1:k}, \mathbf{y}_{k+1:3}))$  is greater than $\pi r^2$,  for all $k=0,1,2$, we have
\begin{equation*}
g_{2,k}(\rho,r,\mathbf{u}_{1:3},\mathbf{y}_{k+1:3}) \leq r^{3}e^{-\frac{1}{2}\pi r^2}\lambda_2(\mathbf{u}_{1:3})\mathbb{1}_{r^2\lambda_2(\mathbf{u}_{1:3})<v_{\rho}}\mathbb{1}_{E_{k,r,\mathbf{u}_{1:3}}}
\end{equation*}
according to \eqref{defunionvol} and \eqref{defg2kdel}. Integrating with respect to $\mathbf{y}_{1:3}$, this gives
\begin{equation}
\label{majG2delminvol}
G_{2,k}(\rho)  \leq \rho\int_{0}^{\infty}\int_{(\mathbf{S}^{1})^{3}}r^{3}e^{-\frac{1}{2}\pi r^2}\lambda_2(\mathbf{u}_{1:3}) \lambda_{2(3-k)}(E_{k,r,\mathbf{u}_{1:3}})\mathbb{1}_{r^2\lambda_2(\mathbf{u}_{1:3})<v_{\rho}}drd\sigma(\mathbf{u}_{1:3}).\end{equation}  
As in the proof of Proposition \ref{maxvoldelappl}, we derive a suitable upper bound of the volume of $E_{k,r,\mathbf{u}_{1:3}}$.

\begin{Le}
\label{volumemaxarea}
Let $\mathbf{u}_{1:3}\in (\mathbf{S}^1)^3$ and $r\geq 0$. Then
\begin{subequations}
\begin{equation}
\label{volumel2}
\lambda_2(E_{2,r,\mathbf{u}_{1:3}}) \leq c\cdot v_{\rho}|u_1-u_2|^{-1}
\end{equation}
\begin{equation}
\label{volumel4}
\lambda_4(E_{1,r,\mathbf{u}_{1:3}}) \leq c\cdot r^2v_{\rho}
\end{equation}
\begin{equation}
\label{volumel6}
\lambda_6(E_{0,r,\mathbf{u}_{1:3}}) \leq c\cdot\log\rho\cdot r^2 v_{\rho}.
\end{equation}
\end{subequations}
\end{Le}

\begin{prooft}{Lemma \ref{volumemaxarea}}
Let $y_3$ be in $E_{2,r,\mathbf{u}_{1:3}}$. Since $R(r\mathbf{u}_{1:2},y_3)$ is less than $r$, we have $|y_3-ru_1|\leq 2R(r\mathbf{u}_{1:2},y_3)\leq 2r$. In particular, we obtain
\begin{equation}
\label{volumeball1}
|y_3|\leq 3r.
\end{equation}
Moreover, the area of the triangle $\Delta(r\mathbf{u}_{1:2},y_3)$ is given by
\begin{equation}\label{areatriangle}\lambda_2(r\mathbf{u}_{1:2},y_3) = \frac{1}{2}r|u_1-u_2|\cdot \delta(y_3,\mathbf{L}(ru_1,ru_2))\end{equation} where $\mathbf{L}(ru_1,ru_2)$ is the affine line induced by the points $p_1=ru_1$, $p_2=ru_2$ and where $\delta(y_3,\mathbf{L}(ru_1,ru_2))$ denotes the distance between this line and the point $y_3$. Since $\lambda_2(r\mathbf{u}_{1:2},y_3)<v_{\rho}$, it results from \eqref{areatriangle} that
\begin{equation}
\label{strip2}
\delta(y_3,\mathbf{L}(ru_1,ru_2))\leq \frac{2v_{\rho}}{r|u_1-u_2|}.
\end{equation}
The inequalities \eqref{volumeball1} and \eqref{strip2} show that $E_{2,r,\mathbf{u}_{1:3}}$ is included in the intersection of a ball of radius $3r$ and a strip of width $\frac{4v_{\rho}}{r|u_1-u_2|}$ i.e. 
\[\lambda_2(E_{2,r,\mathbf{u}_{1:3}}) \leq 6r\times \frac{4v_{\rho}}{r|u_1-u_2|} = c\cdot v_{\rho}|u_1-u_2|^{-1}.\]

Secondly, we bound $\lambda_4(E_{1,r,\mathbf{u}_{1:3}})$. Taking the (spherical coordinates type) change of variables $\phi_2: \RR_+\times\mathbf{S}^1\rightarrow\RR^2, (s',u'_2)\mapsto y_2=ru_1+s'u'_2$ with Jacobian matrix $|D\phi_2 (s',u'_2)| = s'$, we obtain \begin{equation}
\label{l4changeofvariables}
\lambda_4(E_{1,r,\mathbf{u}_{1:3}}) \leq \int_{0}^{2r}\int_{\mathbf{S}^1}\int_{\RR^2}s'\mathbb{1}_{\lambda_2(ru_1,ru_1+s'u'_2, y_3)<v_{\rho}}\mathbb{1}_{R(ru_1,ru_1+s'u'_2, y_3)\leq r}ds'd\sigma(u'_2)dy_3.
\end{equation} The positive number $s'$ is integrated on $[0,2r]$. Indeed, the inequality $R(ru_1,ru_1+s'u'_2, y_3)\leq r$ implies that $s'=|(ru_1+s'u'_2)-ru_1|\leq 2r$. Proceeding along the same lines as in the proof of \eqref{volumel2}, we show that $y_3$ belongs to the ball $B(0,3r)$ and a strip of width $\frac{4v_{\rho}}{s'}$. Integrating \eqref{l4changeofvariables} with respect to $y_3$, we deduce that \[\lambda_4(E_{1,r,\mathbf{u}_{1:3}}) \leq  24\int_{0}^{2r}\int_{\mathbf{S}^1}v_{\rho}rds'd\sigma(u'_2) =  c\cdot r^2 v_{\rho}.\]

Finally, we bound $\lambda_6(E_{0,r,\mathbf{u}_{1:3}})$. Taking the same change of variables as in \eqref{Blaschke}, we have \begin{multline*}\lambda_6(E_{0,r, \mathbf{u}_{1:3}}) \leq \int_{(\RR^2)^3}\mathbb{1}_{z(\mathbf{y}_{1:3})\in \CC}\mathbb{1}_{R(\mathbf{y}_{1:3})<r}\mathbb{1}_{\lambda_2(\mathbf{y}_{1:3})<v_{\rho}}d\mathbf{y}_{1:3}\\ = \int_{\CC}\int_{0}^r\int_{(\mathbf{S}^1)^3}r'^3\lambda_2(\mathbf{u}'_{1:3})\mathbb{1}_{r'^2\lambda_2(u'_1,u'_2,u'_3)<v_{\rho}}dz'dr'd\sigma(\mathbf{u}'_{1:3}).\end{multline*}
Bounding $r'^3\lambda_2(u'_1,u'_2,u'_3)$ by $r'v_{\rho}$ and integrating with respect to $z'\in\CC$, $r'\in [0,r]$ and $\mathbf{u}'_{1:3}\in (\mathbf{S}^1)^3$, we show that $\lambda_6(E_{0,r,\mathbf{u}_{1:3}})$ is less than $c\cdot\lambda_2(\CC)r^2v_{\rho}$ with $\lambda_2(\CC)\leq c\cdot\log\rho$. 
\end{prooft}

We can now derive a suitable upper bound of $G_{2,k}(\rho)$. Indeed, if $k=0$, we deduce from \eqref{majG2delminvol} and \eqref{volumel6} that 
\begin{multline*}G_{2,0}(\rho) \leq c\cdot\log\rho\cdot \rho v_{\rho}\int_{0}^{\infty}\int_{(\mathbf{S}^1)^3}r^{5}e^{-\frac{1}{2}\pi r^2}\lambda_2(\mathbf{u}_{1:3}) \mathbb{1}_{r^2\lambda_2(\mathbf{u}_{1:3})<v_{\rho}}drd\sigma(\mathbf{u}_{1:3})\\
\leq c\cdot\log\rho\cdot \rho v_{\rho}^2\int_{0}^{\infty}\int_{(\mathbf{S}^1)^3}r^{3}e^{-\frac{1}{2}\pi r^2}drd\sigma(\mathbf{u}_{1:3}).\end{multline*} First, we notice that the integral of the right-hand side is bounded. Moreover, replacing $v_{\rho}$ by $c\cdot \rho^{-3/5}$ according to \eqref{defvminvol}, we show that $G_{2,0}(\rho)$ is less than $c\cdot \log\rho\cdot \rho^{-1/5}$. In the same spirit, when $k=1$, we obtain that $G_{2,1}(\rho)\leq c\cdot \rho^{-1/5}$ according to \eqref{majG2delminvol} and \eqref{volumel4}. Hence
\begin{equation}
\label{majG20G21minvol}
G_{2,0}(\rho) = O\left(\log\rho\cdot \rho^{-1/5}\right) \text{ and } G_{2,1}(\rho) = O\left(\rho^{-1/5}\right).
\end{equation}

Finally, if $k=2$, we deduce from \eqref{majG2delminvol} and \eqref{volumel2} that
\begin{equation*}
G_{2,2}(\rho) \leq c\cdot\rho v_{\rho}\int_{0}^{\infty}\int_{(\mathbf{S}^1)^3}r^3e^{-\frac{1}{2}\pi r^2}\lambda_2(\mathbf{u}_{1:3})|u_1-u_2|^{-1}\mathbb{1}_{r^2\lambda_2(\mathbf{u}_{1:3})<v_{\rho}}drd\sigma(\mathbf{u}_{1:3}).
\end{equation*} Let $\phi_3$ be the change of variables 
\begin{equation*} \begin{split} \phi_3: & [0,2\pi)^3\longrightarrow (\mathbf{S}^1)^3\\
& \boldsymbol{\theta}_{1:3} \longmapsto \mathbf{u}_{1:3} \text{ with } u_1 = u(-\theta_1+\theta_3),  u_2=u(\theta_1+\theta_3) \text{ and } u_3= u(\theta_2+\theta_3) \end{split}\end{equation*} where $u(\theta) = (\cos\theta, \sin\theta)$. For all $\boldsymbol{\theta}_{1:3}\in [0,2\pi)^3$, let us denote by $A(\boldsymbol{\theta}_{1:3}) = \lambda_2(\mathbf{u}_{1:3})$ with $\mathbf{u}_{1:3} = \phi_3(\boldsymbol{\theta}_{1:3})$. Since $|u_1-u_2| = 2|\sin\theta_1|$, we have
\begin{equation*}
G_{2,2}(\rho) \leq c\cdot\rho v_{\rho}\int_{0}^{\infty}\int_{[0,\pi/2)\times [0,2\pi)^2}r^3e^{-\frac{1}{2}\pi r^2}A(\boldsymbol{\theta}_{1:3})|\sin\theta_1|^{-1}\mathbb{1}_{r^2A(\boldsymbol{\theta}_{1:3})<v_{\rho}}drd\boldsymbol{\theta}_{1:3}.
\end{equation*} 
Without loss of generality, we have assumed that $\theta_1$ belongs to $[0,\pi/2]$. To bound $G_{2,2}(\rho)$, we consider two cases that depend on the order of $\theta_1$. Let $\epsilon>\frac{3}{5}$ be fixed. The previous inequality can be written as 
\begin{multline}
\label{majG22minvol}
G_{2,2}(\rho) \leq c\cdot\rho v_{\rho}\int_{0}^{\infty}\int_{[0,\rho^{-\epsilon}[\times [0,2\pi)^2}r^3e^{-\frac{1}{2}\pi r^2}A(\boldsymbol{\theta}_{1:3})|\sin\theta_1|^{-1}\mathbb{1}_{r^2A(\boldsymbol{\theta}_{1:3})<v_{\rho}}drd\boldsymbol{\theta}_{1:3}\\
+  c\cdot\rho v_{\rho}\int_{0}^{\infty}\int_{[\rho^{-\epsilon},\pi/2)\times [0,2\pi)^2}r^3e^{-\frac{1}{2}\pi r^2}A(\boldsymbol{\theta}_{1:3})|\sin\theta_1|^{-1}\mathbb{1}_{r^2A(\boldsymbol{\theta}_{1:3})<v_{\rho}}drd\boldsymbol{\theta}_{1:3} = G_{2,2}^{(1)}(\rho)  + G_{2,2}^{(2)}(\rho)
\end{multline} where $G_{2,2}^{(1)}(\rho)$ and $G_{2,2}^{(2)}(\rho)$ denote respectively the first and the second integrals of the right-hand side. Let us note that $A(\boldsymbol{\theta}_{1:3})|\sin\theta_1|^{-1}$ is bounded since, according to \eqref{areatriangle}, we have $A(\boldsymbol{\theta}_{1:3}) = \frac{1}{2}\cdot 2|\sin\theta_1|\cdot d(u_3,\mathbf{L}(\mathbf{u}_{1:2}))$ where $\mathbf{u}_{1:3} = \phi_3(\boldsymbol{\theta}_{1:3})$ and $d(u_3,\mathbf{L}(\mathbf{u}_{1:2}))\leq 2$. Hence, the first integral of the right-hand side of \eqref{majG22minvol} is less than
\begin{equation}
\label{majG22minvolint1}
G_{2,2}^{(1)}(\rho) \leq c\cdot \rho v_{\rho}\int_{0}^{\infty}\int_{[0,\rho^{-\epsilon})\times [0,2\pi)^2}r^3e^{-\frac{1}{2}\pi r^2}drd\boldsymbol{\theta}_{1:3}\leq c\cdot\rho^{1-\epsilon}v_{\rho} = O(\rho^{-1/5})
\end{equation} since $v_{\rho}=c\cdot\rho^{-3/5}$ and $\epsilon>\frac{3}{5}$. Moreover, bounding $A(\boldsymbol{\theta}_{1:3})$ by $r^{-2}v_{\rho}$ in the second integral of \eqref{majG22minvol}, we have
\begin{equation}
\label{majG22minvolint2}
G_{2,2}^{(2)}(\rho) \leq  c\cdot\rho v_{\rho}^2\int_{0}^{\infty}\int_{[\rho^{-\epsilon},\pi/2)\times [0,2\pi)^2}re^{-\frac{1}{2}\pi r^2}|\sin\theta_1|^{-1}drd\boldsymbol{\theta}_{1:3}\leq c\cdot \log\rho\cdot\rho v_{\rho}^{2} = O\left(\log\rho\cdot \rho^{-1/5} \right)
\end{equation} since $\int_{\rho^{-\epsilon}}^{\pi/2}\frac{1}{|\sin\theta_1|}d\theta_1$ is of order $\log\rho$. 

From  \eqref{majG20G21minvol}, \eqref{majG22minvol}, \eqref{majG22minvolint1} and \eqref{majG22minvolint2}, we deduce that $G_2(\rho) = O\left(\log\rho\cdot\rho^{-1/5}\right)$. Proposition \ref{minvoldelappl} is now a direct consequence of \eqref{majG1minvol} and Theorem \ref{realtyp}.
\end{prooft}

The main tool to derive the asymptotic behaviour of $A_{PDT,\min}(\rho)$ is the Taylor expansion of $K_{1/6}(\cdot)$ used in \eqref{Taylorminvol}. To the best of our knowledge, there is not more accurate result on this Taylor expansion which could provide the rate of convergence  $\PPP{\lambda_2(\cell)<v}$. Actually, the rate of convergence can be investigated with a more complicated method. Indeed, in \cite{Ra}, using Mellin transform, Rathie shows that the density of $\lambda_2(\cell)$ is given by
\[f(x) = 3\pi^{-1/2}(2\pi ix)^{-1}\int_{L}\frac{\Gamma(z+5/6)\Gamma(z+1)\Gamma(z+7/6)}{\Gamma(z+3/2)}(4\pi x^2/27)^{-z}dz\] where $L$ encloses all the (complex) poles of the integrand. These poles, of order 1, are $-5/6-k$, $-1-k$ and $-7/6-k$, $k=0,1,2,\cdots$. Evaluating the contour integral as the sum of the residues at the poles, he shows that
\[f(x) = \sum_{k=0}^{\infty}c_{k,1}x^{2/3+2k} + \sum_{k=0}^{\infty}c_{k,2}x^{1+2k} + \sum_{k=0}^{\infty}c_{k,3}x^{4/3+2k}.\] It results of a Taylor expansion of the sums that $f(x) = c_{0,1}x^{2/3}+O(x)$. Integrating $f(\cdot)$ on $[0,v]$, we obtain that 
\[\PPP{\lambda_2(\cell)<v} = c_{0,1}\cdot v^{5/3}+O(v^2).\] Taking $v=v_{\rho}$ as in \eqref{defvminvol}, the function
$G_1(\rho) = |\rho\PPP{\lambda_2(\cell)<v_{\rho}} - t^{5/3}|$ is of order $\rho v_{\rho}^2 = c\cdot\rho^{-1/5}$. Since $G_2(\rho) = O\left(\log\rho\cdot\rho^{-1/5}\right)$, we obtain the more precise result
\begin{equation*} \left|\PPP{\alpha_{3}^{3/5}\rho^{3/5}A_{PDT,\min}(\rho)\geq t} - e^{-t^{5/3}}\right| = O\left(\log\rho\cdot\rho^{-1/5} \right). \end{equation*} Nevertheless, we have used the Taylor expansion of the modified Bessel function to prove Proposition \ref{minvoldelappl} since the method is quicker than the use of series. 

When $d\geq3$, the density of $\lambda_3(\cell)$ can also be written as an integral (see (2.5) in \cite{Ra}):
\[f(x) = c_1(2\pi ix)^{-1}\int_{L}\bigtriangledown_d(z)\cdot(c_2x^2)^{-z}dz\] where $c_1$, $c_2$ are two constants depending on $d$ which can be made explicit and
\[\bigtriangledown_d(z) = \frac{\prod_{j=2}^d\Gamma(j/2+z)\prod_{j=0}^d\Gamma\left(\frac{d^2+1+2j}{2(d+1)}+z \right)    }{\prod_{j=1}^{d-1}\Gamma(d/2+j/d+z)\Gamma^{d-1}((d+1)/2+z)}.\] The poles of $\bigtriangledown_d(\cdot)$ are real numbers and the largest of them is $-1$ which is a simple pole. Proceeding along the same lines as in the case $d=2$, we show that $f(x) = c\cdot x+o(x)$ when $x$ goes to 0 i.e. 
\[G_1(\rho) = \left|\rho\PPP{\lambda_d(\cell)<c\cdot\rho^{-1/2}t} - t^2\right| \conv[\rho]{\infty}0\] for $d\geq 3$. Unfortunately, the same method as in the proof of Proposition \ref{minvoldelappl} is not enough to show that $G_2(\rho)$ converges to 0. Nevertheless, we would be able to show that the minimum of the volumes of the cells of a Poisson-Delaunay tessellation is of order $\rho^{-1/2}$ provided that the extremal index exists and differs from 0 (see section \ref{sectionrealtyp2} for more details about extremal index).

\section{Extreme Values of a Poisson-Voronoi tessellation}
\label{sectionPVT}
Let $\chi$ be a locally finite subset of $\RR^d$. For all $x\in\chi$, we denote by $C_{\chi}(x)$ the Voronoi cell of nucleus $x$ defined as 
\[C_{\chi}(x) = \{y\in\RR^d, |x-y|\leq |x'-y|, x'\neq x\in \chi\}.\] For all $x\in\chi$, we denote by $\mathcal{N}_{\chi}(x)$ the set of neighbors of $x$ and $N_{\chi}(x)$ its cardinality i.e. \begin{equation}\label{notneighbor}\mathcal{N}_{\chi}(x) = \{x'\in \mathbf{X}, C_{\chi}(x')\cap C_{\chi}(x)\neq\varnothing\} \text { and } N_{\chi}(x)=\#\mathcal{N}_{\chi}(x).\end{equation} 
Voronoi tessellation corresponds to the dual graph of Delaunay tessellation in the following sense: there exists an edge between two points $x,x'\in\chi$ in the Delaunay graph if and only if they are Voronoi neighbors i.e. $C_{\chi}(x)\cap C_{\chi}(x')\neq\varnothing$.

When $\chi = \mathbf{X}$ is a Poisson point process (of intensity 1), the family $\mm_{PVT}=\{C_{\mathbf{X}}(x), x\in \mathbf{X}\}$ is called the Poisson-Voronoi tessellation. Such model is extensively used in many domains such as cellular biology \cite{Po}, astrophysics \cite{RBFN}, telecommunications \cite{BB3} and ecology \cite{Ro}. For a complete account, we refer to the books \cite{M2}, \cite{OBSC}, \cite{SW}  and the survey \cite{Cal5}. 

As in section \ref{sectionPVT}, the window $\mathbf{W}_{\rho}=\rho^{1/d}[0,1]^d$ is partitioned into $N_{\rho} = \left\lfloor\frac{\rho}{2\log\rho}\right\rfloor$ sub-cubes $\mathbf{i}\in V_{\rho}$. The event $A_{\rho}$ is the same as in \eqref{defAgammadel} and we can show that it satisfies \textsc{Condition 1} for the Poisson-Voronoi tessellation with  arguments very similar to the proof of Lemma \ref{AssumptAgammadel}. 

For each cell $C\in\mm_{PVT}$ i.e. $C=C_{\mathbf{X}}(x)$, we take  $z(C_{\mathbf{X}}(x))=x$. A consequence of Slivnyak's Theorem (see e.g. Theorem 3.3.5  in \cite{SW}) shows that the typical cell satisfies the equality in distribution \begin{equation}\label{typicalcellvoronoi}\cell \overset{\mathcal{D}}{=} C_{\mathbf{X}\cup\{0\}}(0)\end{equation} where   $C_{\mathbf{X}\cup\{0\}}(0)$ is the Voronoi cell of nucleus $0$ when we add the origin to the Poisson point process.  

The function $G_2(\cdot)$ defined in \eqref{defG2general} has an integral representation. Indeed, from Slivnyak's Formula, it can be written as \begin{equation}
\label{genericVoronoi}
G_2(\rho) = \rho\int_{\CC}\PPP{f(C_{\mathbf{X}\cup\{0,y\}}(0))>v_{\rho}, f(C_{\mathbf{X}\cup\{0,y\}}(y))>v_{\rho}}dy.
\end{equation}

Extremes of characteristic radii of Poisson-Voronoi tessellation are studied in \cite{CC}. In this paper, we give the asymptotic behaviours of two new geometrical characteristics. 

The first one is the distance to the farthest neighbor. More precisely, we consider \begin{equation}\label{deffarthest}D(C_{\mathbf{X}}(x)) = \max_{x'\in \mathcal{N}_{\mathbf{X}}(x)}|x-x'|, x\in \mathbf{X} \text{ and } D_{\min,PVT}(\rho) = \min_{x\in \mathbf{X}\cap \mathbf{W}_{\rho}}D(C_{\mathbf{X}}(x)).\end{equation}
The second characteristic is the volume of the so-called Voronoi flower. We denote respectively for each point $x\in \mathbf{X}$, the Voronoi flower of nucleus $x$ and the minimum of their volumes  as \begin{equation}\label{defflower}\mathcal{F}(C_{\mathbf{X}}(x)) = \bigcup_{y\in C_{\mathbf{X}}(x)}B(y, |y-x|) \text{ and } F_{\min,PVT}(\rho) = \min_{x\in\mathbf{X}\cap \mathbf{W}_{\rho}}\lambda_d(\mathcal{F}(C_{\mathbf{X}}(x))).\end{equation}   Obviously, $2^{-d}\kappa_d D^d_{\min,PVT}(\rho)\leq\kappa_d \min_{x\in\mathbf{X}\cap \mathbf{W}_{\rho}}R(C_{\mathbf{X}}(x))^d\leq F_{\min,PVT}(\rho)$ where $R(C_{\mathbf{X}}(x))$ denotes the circumradius of $C_{\mathbf{X}}(x)$. Actually, the following proposition shows that the two random variables  $D^d_{\min,PVT}(\rho)$ and $F_{\min,PVT}(\rho)$   are of same order when $\rho$ goes to infinity.

\begin{Prop}
\label{VoronoiProp}
Let $\mm_{PVT}$ be a Poisson-Voronoi tessellation of intensity $\gamma=1$. For all $t\geq 0$, we have
\begin{subequations}
      \begin{equation}\label{farthestneighbor}
  \left|\PPP{\alpha_{d,4}^{1/(d+1)}\rho^{1/(d+1)}D^d_{\min,PVT}(\rho) \geq t  } - e^{-t^{d+1}}\right| = O\left(\rho^{-1/(d+1)}\right) \end{equation}

  \begin{equation}\label{flower}
  \left|\PPP{\alpha_{d,5}^{1/(d+1)}\rho^{1/(d+1)}F_{\min,PVT}(\rho) \geq t  } - e^{-t^{d+1}}\right| = O\left(\rho^{-1/(d+1)}\right) \end{equation}
    \end{subequations}
 where $\alpha_{d,4}$ and $\alpha_{d,5}$ are given in \eqref{defalpha1} and \eqref{defalpha2} respectively.
 
\end{Prop}

Before proving Proposition \ref{VoronoiProp}, we need a practical lemma which is a new version of Lemma 3 in \cite{CC} adapted to our framework. 
\begin{Le}
\label{boundedcelllemma}
Let $v\geq 0$, $y\neq 0\in\RR^d$ and $\chi\subset\RR^d$ locally finite such that $\chi\cup\{0,y\}$ is in general position i.e. each subset of size n<$d+1$ is affinely independent (see \cite{Ze}). Let us assume that each Voronoi cell associated to the set $\chi\cup\{0,y\}$ is bounded and that \begin{equation}\label{condboundedcelllemma} \mathcal{N}_{\chi\cup\{0,y\}}(0)\subset B(0,v)\text{ and } \mathcal{N}_{\chi\cup\{0,y\}}(y)\subset B(y,v).\end{equation} Then
\[\#\left(\chi\cap \left(B(0,v)\cup B(y,v) \right)\right)\geq d+1.\]
\end{Le}

\begin{prooft}{Lemma \ref{boundedcelllemma}}
Let us define $\chi_{0,y}$ as the (finite) subset: 
\[\chi_{0,y} = \chi \cap \left(B(0,v)\cup B(y,v)\right).\] Thanks to \eqref{condboundedcelllemma}, we have $C_{\chi\cup\{0,y\}}(0) = C_{\chi_{0,y}\cup \{0,y\}}(0)$ and $C_{\chi\cup\{0,y\}}(y) = C_{\chi_{0,y}\cup \{0,y\}}(y)$. In  particular, this shows that the cells  $C_{\chi_{0,y}\cup \{0,y\}}(0)$ and $C_{\chi_{0,y}\cup \{0,y\}}(y)$ are bounded. Hence $0$ and $y$ are in the convex hulls of $\chi_{0,y}\cup\{y\}$ and $\chi_{0,y}\cup\{0\}$ respectively (see Property \textbf{V}2, p. 58 in \cite{OBSC}). This implies that \[\{0,y\}\subset \text{conv}(\chi_{0,y}).\] Since $\chi\cup\{0,y\}$ is in general position, this shows that $\text{conv}(\chi_{0,y})$ has a non-empty interior and consequently this proves Lemma \ref{boundedcelllemma}.
 \end{prooft}

We can now prove Proposition \ref{VoronoiProp}.

\begin{prooft}{Proposition \ref{VoronoiProp}}

\textbf{Proof of \eqref{farthestneighbor}.} To find a function $v_{\rho}(t)$ such that $G_1(\rho) = | \rho\PPP{D(\cell)>v_{\rho}} -t |$ converges to 0, we have to approximate the tail of the distribution function of $D(\cell)$. Let $v\geq 0$ be fixed. Since $\cell = C_{\mathbf{X}\cup\{0\}}(0)$, we have \begin{equation}\label{inclfarthest}D(\cell)<v \Longleftrightarrow \mathcal{N}_{\mathbf{X}\cup\{0\}}(0)\subset B(0,v).\end{equation} In particular, we get
 \begin{equation}\label{Conditioningfarthest}\PPP{D(\cell)<v} = \sum_{k=d+1}^{\infty}\PPP{\mathcal{N}_{\mathbf{X}\cup\{0\}}(0)\subset B(0,v), N_{\mathbf{X}\cup\{0\}}(0) = k}.\end{equation}  An integral representation of the right-hand side is given by (see Proposition 1 in \cite{Cal7}) \begin{equation*}
  \PPP{\mathcal{N}_{\mathbf{X}\cup\{0\}}(0)\subset B(0,v),N_{\mathbf{X}\cup\{0\}}(0) = k} 
  = \frac{1}{k!}\int_{B(0,v)^k}e^{-\lambda_d(\mathcal{F}(C_{\{\mathbf{x}_{1:k}\}\cup\{0\}}(0)))} \mathbb{1}_{F_k}(\mathbf{x}_{1:k})d\mathbf{x}_{1:k}
 \end{equation*} where \[F_k = \left\{\mathbf{x}_{1:k}=(x_1,\ldots, x_k)\in (\RR^d)^k, C_{\{\mathbf{x}_{1:k}\} \cup\{0\}}(0) \text{ is a convex polytope with $k$ faces}\right\}.\] We recall that $\{\mathbf{x}_{1:k}\} \cup\{0\} = \{x_1, x_2, \ldots, x_k, 0\}$.  Taking the change of variables $x_i=vx'_i$, we obtain for all $k\geq d+1$ \begin{equation}\label{probafarthest}\PPP{\mathcal{N}_{\mathbf{X}\cup\{0\}}(0)\subset B(0,v),N_{\mathbf{X}\cup\{0\}}(0) = k} = v^{dk}\cdot\frac{1}{k!}\int_{B(0,1)^{k}}e^{-v^d\lambda_d(\mathcal{F}(C_{\{\mathbf{x}'_{1:k}\}\cup\{0\}}(0)))}\mathbb{1}_{F_k}(\mathbf{x}'_{1:k})d\mathbf{x}'_{1:k}.\end{equation} If $k=d+1$, the previous probability converges to $\alpha_{d,4}\cdot v^{d(d+1)}$ when $v$ goes to 0 where \begin{equation}
 \label{defalpha1}
 \alpha_{d,4} = \frac{1}{(d+1)!}\int_{B(0,1)^{d+1}}\mathbb{1}_{A_{d+1}}(\mathbf{x}'_{1:d+1})d\mathbf{x}'_{1:d+1}.
 \end{equation} If $k\geq d+2$, the right-hand side of \eqref{probafarthest} is less than $\frac{\kappa_d^k}{k!}v^{dk}$ thanks to the trivial inequalities $\mathbb{1}_{F_k}\leq 1$ and $e^{-\lambda_d(\mathcal{F}(C_{\{\mathbf{x}'_{1:k}\}\cup\{0\}}(0)))}\leq 1$. It follows from \eqref{Conditioningfarthest} that 
 \begin{equation}\label{taylorfarthest}\left|\PPP{D(\cell)<v} - \alpha_{d,4}\cdot v^{d(d+1)} \right| \leq \sum_{k=d+2}^{\infty}\frac{\kappa_d^k}{k!}v^{dk} = O(v^{d(d+2)}).\end{equation} Now, we can choose a suitable function $v_{\rho}$. Indeed, let $t \geq 0$ be fixed and let us denote by 
 \begin{equation}
 \label{defvfarthest}
 v_{\rho} = v_{\rho}(t) = \left(\alpha_{d,4}^{-1}\rho^{-1} \right)^{1/d(d+1)}t^{1/d}.
 \end{equation} According to \eqref{taylorfarthest}, we have \begin{equation}\label{majG1farthest}G_1(\rho) = |\rho\PPP{D(\cell)<v_{\rho}} - t^{d+1} | = O\left(\rho^{-1/(d+1)}\right).\end{equation} 
 
Let us give now an upper bound of the function $G_2(\rho)$ defined in \eqref{defG2general}. According to  \eqref{genericVoronoi} and in the same spirit as in \eqref{inclfarthest}, we obtain that
\begin{multline}
\label{majG2farthest0}
G_2(\rho) = \rho\int_{\CC}\PPP{D(C_{\mathbf{X}\cup\{0,y\}}(0))<v_{\rho},D(C_{\mathbf{X}\cup\{0,y\}}(y))<v_{\rho} }dy\\= \rho\int_{\CC}\PPP{\mathcal{N}_{\mathbf{X}\cup\{0,y\}}(0)\subset B(0,v_{\rho}), \mathcal{N}_{\mathbf{X}\cup\{0,y\}}(y)\subset B(y,v_{\rho})}dy. \end{multline} 
To guarantee the independence of the events considered in \eqref{majG2farthest0} for each cells which are distant enough, we write
\begin{multline}
\label{majG2farthest0bis}
G_2(\rho) = \rho\int_{\CC\cap B(0,2v_{\rho})^c}\PPP{\mathcal{N}_{\mathbf{X}\cup\{0,y\}}(0)\subset B(0,v_{\rho}), \mathcal{N}_{\mathbf{X}\cup\{0,y\}}(y)\subset B(y,v_{\rho})}dy\\
+ \rho\int_{\CC\cap B(0,2v_{\rho})}\PPP{\mathcal{N}_{\mathbf{X}\cup\{0,y\}}(0)\subset B(0,v_{\rho}), \mathcal{N}_{\mathbf{X}\cup\{0,y\}}(y)\subset B(y,v_{\rho})}dy. 
\end{multline}

For the first integral, when $y\in \CC\cap B(0,2v_{\rho})^c$, the balls $B(0,v_{\rho})$ and $B(y, v_{\rho})$ are disjoint. Because $\mathbf{X}$ is a Poisson point process and because $y\not\in B(0,2v_{\rho})$, the first integrand of \eqref{majG2farthest0bis} can be written as the product $\PPP{\mathcal{N}_{\mathbf{X}\cup\{0\}}(0)\subset B(0,v_{\rho})}\times \PPP{ \mathcal{N}_{\mathbf{X}\cup\{y\}}(y)\subset B(y,v_{\rho})}$. Hence, according to \eqref{inclfarthest} and \eqref{majG1farthest} we obtain that
\begin{equation}\label{integrandfarthest}\PPP{\mathcal{N}_{\mathbf{X}\cup\{0,y\}}(0)\subset B(0,v_{\rho}), \mathcal{N}_{\mathbf{X}\cup\{0,y\}}(y)\subset B(y,v_{\rho})} = \PPP{D(\cell)<v_{\rho}}^2\leq c\cdot \rho^{-2}, \hspace{0.2cm} y\in B(0,2v_{\rho})^c\end{equation} where $c$ is a constant which \textit{does not} depend on $y$.

For the second integral of \eqref{majG2farthest0bis}, we apply Lemma \ref{boundedcelllemma} to $\chi=\mathbf{X}$. This gives \begin{equation}\label{poissfarthest1}\PPP{\mathcal{N}_{\mathbf{X}\cup\{0,y\}}(0)\subset B(0,v_{\rho}), \mathcal{N}_{\mathbf{X}\cup\{0,y\}}(y)\subset B(y,v_{\rho})}\leq \PPP{\#(\mathbf{X}\cap (B(0,v_{\rho})\cup B(y,v_{\rho}))\geq d+1)}, \hspace{0.2cm} y\in B(0,2v_{\rho}).\end{equation} Since $\#(\mathbf{X}\cap B)$ is Poisson distributed of mean $\lambda_d(B)$ for each Borel subset $B\subset\RR^d$, we obtain  for $\rho$ large enough that 
\begin{multline*}\label{poissfarthest2}\PPP{\#(\mathbf{X}\cap (B(0,v_{\rho})\cup B(y,v_{\rho}))\geq d+1)} = \sum_{k=d+1}^{\infty}\frac{1}{k!}\left(\lambda_d(B(0,v_{\rho})\cup B(y,v_{\rho}))\right)^ke^{-\lambda_d(B(0,v_{\rho})\cup B(y,v_{\rho}))}\\\leq c\cdot v_{\rho}^{d(d+1)} = c'\cdot \rho^{-1}, y\in B(0,2v_{\rho})\end{multline*} according to \eqref{defvfarthest} and to the trivial inequalities  $e^{-\lambda_d(B(0,v_{\rho})\cup B(y,v_{\rho}))}\leq 1$ and $\lambda_d(B(0,v_{\rho})\cup B(y,v_{\rho}))\leq 2\cdot \kappa_d v^d_{\rho}$. This together with \eqref{majG2farthest0bis}, \eqref{integrandfarthest} and \eqref{poissfarthest1} shows that
\[G_2(\rho) \leq c\cdot \rho^{-1}\lambda_d(\mathfrak{C}_{\rho}\cap B(0,2v_{\rho})^c) + c\cdot \lambda_d(\mathfrak{C}_{\rho}\cap B(0,2v_{\rho})).\] Since $\lambda_d(\CC\cap B(0,2v_{\rho})^c)\leq \lambda_d(\CC)\leq c\cdot\log\rho$ and $\lambda_d(\CC\cap B(0,2v_{\rho}))\leq \lambda_d(B(0,2v_{\rho})) = c\cdot \rho^{-1/(d+1)}$, we deduce from the previous inequality that \begin{equation}\label{majG2farthest}G_2(\rho) \leq c\cdot\log\rho\times\rho^{-1} + c\cdot \rho^{-1/(d+1)} = O\left(\rho^{-1/(d+1)} \right).\end{equation} We now derive directly \eqref{farthestneighbor} from \eqref{majG1farthest}, \eqref{majG2farthest} and Theorem \ref{realtyp}. 

\vspace{0.5cm}

\textbf{Proof of \eqref{flower}}. This will be sketched since it is analogous to the proof of \eqref{farthestneighbor}. First, we investigate the tail of the distribution function of $\lambda_d(\mathcal{F}(\cell))$. In \cite{Z}, Zuyev shows that, conditional on $N_{\mathbf{X}\cup\{0\}}=k$, the volume of $\mathcal{F}(\cell)$ is Gamma distributed of parameters $(k,1)$ i.e. 
\begin{equation}
\label{Conditioningflower}
 \PPP{\lambda_d(\mathcal{F}(\cell))< v} = \sum_{k=d+1}^{\infty}\frac{1}{(k-1)!}\int_0^{ v}x^{k-1}e^{-x}dx\cdot p(k)  
\end{equation} where $p(k)=\PPP{N_{\mathbf{X}\cup\{0\}}(0)=k}$. When $k=d+1$, the Taylor expansion $e^{-x} = 1+O(x)$ shows that the term of the series in \eqref{Conditioningflower} equals $\alpha_{d,5}v^{d+1} + O(v^{d+2})$ where \begin{equation}
\label{defalpha2}
\alpha_{d,5} = \frac{p(d+1)}{(d+1)!}.
\end{equation} If $k\geq d+2$, the term of the series in \eqref{Conditioningflower} is less than $\frac{1}{d!}\cdot v^{d+2}\cdot p(k)$ thanks to the trivial inequality $e^{-x}\leq 1$. According to \eqref{Conditioningflower}, we get \[| \PPP{\lambda_d(\mathcal{F}(\cell))< v} - \alpha_{d,5}\cdot v^{d+1}| = O(v^{d+2}).\] Hence, for all fixed $t\geq 0$, taking \begin{equation}\label{defvflower}v_{\rho} = v_{\rho}(t) = \left(\alpha_{d,5}^{-1}\rho^{-1} \right)^{1/(d+1)}t\end{equation} we obtain \begin{equation} \label{majG1flower}G_1(\rho) = | \rho\PPP{\lambda_d(\mathcal{F}(\cell))< v_{\rho}} - t^{d+1}| = O(\rho^{-1/(d+1)}).\end{equation}

To get an upper bound of $G_2(\rho)$, we note that for each $\chi\subset \RR^d$ locally finite and $x\in\chi$, we have
\[\frac{\kappa_d}{2^d}\cdot \left(D(C_{\chi}(x)) \right)^d \leq \lambda_d\left(\mathcal{F}(C_{\chi}(x)) \right)\] where $D(C_{\chi}(x))$ and $\mathcal{F}(C_{\chi}(x))$ are defined as in \eqref{deffarthest} and \eqref{defflower}. Applying the previous inequality to $\chi = \mathbf{X}\cup\{0,y\}$ and $x=0,y$, we deduce from \eqref{genericVoronoi} that

\begin{multline}
\label{majG2flower1}
G_2(\rho) =  \rho\int_{\CC}\PPP{\lambda_d(\mathcal{F}(C_{\mathbf{X}\cup\{0,y\}}(0)))<v_{\rho},\lambda_d(\mathcal{F}(C_{\mathbf{X}\cup\{0,y\}}(y)))<v_{\rho} }dy\\
\leq \rho\int_{\CC}\PPP{D(C_{\mathbf{X}\cup\{0,y\}}(0))<v'_{\rho},D(C_{\mathbf{X}\cup\{0,y\}}(y))<v'_{\rho} }dy
\end{multline} with \[v'_{\rho} = 2\kappa_d^{1/d}\cdot v_{\rho}^{1/d} = (2^{d(d+1)}\kappa_d^{d+1}\alpha_{d,5}^{-1}\rho^{-1})^{1/d(d+1)}t^{1/d}\] according to \eqref{defvflower}.
Let us notice that the function $v'_{\rho}$ equals $v_{\rho}$, defined in \eqref{defvfarthest}, up to a multiplicative constant. Writing the right-hand side of \eqref{majG2flower1} in the same spirit as in \eqref{majG2farthest0} and proceeding along the same lines as in the proof of \eqref{farthestneighbor}, we show that $G_2(\rho)$ is of order $\rho^{-1/(d+1)}$. This together with \eqref{majG1flower} shows \eqref{flower}.

\end{prooft}

The random variables $F_{\min,PVT}(\rho)$ and $D_{\min,PVT}(\rho)$ are related to the minimum of the circumradii $R_{\min, PVT}(\rho)$ which is defined in \cite{CC} since both investigate a minimax. In the same spirit as before, we could re-find the asymptotic behaviour of $R_{\min, PVT}(\rho)$ included in \cite{CC} and prove that the rate of convergence is of order $\rho^{-1/(d+1)}$.

\section{The maximum of inradii of a Gauss-Poisson Voronoi tessellation}
\label{sectionGP}
As an example of non-Poisson point process, a Gauss-Poisson process is analyzed. Introduced by Newman and investigated by Milne and Westcott, such process has a potential application in statistical mechanics (see \cite{Ne}, p. 350) and could be used as a model for molecular motion (see \cite{MWe} p. 169). In the sense of \cite{SKM} p. 161, a stationary planar Gauss-Poisson process $\mathbf{X}$ is a (simple) point process which can be defined as follows: let $\mathbf{X}_a$ be a Poisson point process of intensity $\gamma_a$ in $\RR^2$. Every point $x_a\in\mathbf{X}_a$ is replaced by a cluster of points $\Xi(x_a)=x_a+\Xi_0(x_a)$  where the set of points $\Xi_0(x_a), x_a\in\mathbf{X}_a$ are chosen independently and with identical distribution i.e. 
\[\mathbf{X} = \bigcup_{x_a\in \mathbf{X}_a}\Xi(x_a).\] For all $x_a\in\mathbf{X}_a$, the cluster $\Xi_0(x_a)$ equals in distribution $\Xi_0$ which is defined in the following sense: $\Xi_0$ has an isotropic distribution and is composed of zero, one or two points with probability $p_0\neq 1$, $p_1$ and $p_2=1-(p_0+p_1)$. If $\Xi_0$ contains only one point then that point is the origin 0. If $\Xi_0$ is composed of two points then these are separated by a unit distance and have midpoint 0. The intensity of $\mathbf{X}$ is given by 
\[\gamma_{\mathbf{X}} = (p_1+2p_2)\cdot \gamma_a.\] In this subsection, we investigate the maximum of inradii of a Gauss-Poisson Voronoi tessellation $\mm_{GPVT}$ i.e.
\[r_{\max, GPVT}(\rho) = \max_{x\in \mathbf{X}\cap \mathbf{W}_{\rho}}r(C_{\mathbf{X}}(x)) \hspace*{0.2cm}\text{ where } \hspace{0.2cm}
r(C_{\mathbf{X}}(x)) = \max\{r\geq 0, B(x,r)\subset C_{\mathbf{X}}(x)\}.\]
To apply Theorem \ref{realtyp}, we subdivide $\mathbf{W}_{\rho}$ into $N_{\rho}$ sub-cubes of equal size where we take \[N_{\rho} = \left\lfloor \frac{\gamma_a(p_1+p_2) \rho}{2\log\rho}  \right\rfloor.\] With the same method as for a Poisson-Voronoi tessellation, we can show that there exists an integer $R\geq 1$ and a event $A_{\rho}$ (in the same spirit as in \eqref{defAgammadel}) such that \textsc{Condition 1} holds when the Voronoi tessellation is induced by a Gauss-Poisson process. The asymptotic distribution of $r_{\max}(\rho)$  is given in the following proposition. 
\begin{Prop}
\label{GaussPoissonprop}
Let $\mathbf{X}$ be a Gauss-Poisson process of intensity 1 i.e. $(p_1+2p_2)\gamma_a=1$ with $p_0\neq 1$ and $p_1\neq 0$. For all $t\in\RR$, we have 
\[\left|\PPP{r_{\max, GPVT}(\rho)\leq v_{\rho}} - e^{-e^{-t}}\right| = O\left((\log\rho)^{-1/2}\right)\] where $v_{\rho}=v_{\rho}(t)$ is given in \eqref{defvGP}.
\end{Prop}

\begin{prooft}{Proposition \ref{GaussPoissonprop}}
We notice that for all $x\in\mathbf{X}$ and $v\geq 0$, the inscribed radius $r(C_{\mathbf{X}}(x))$ is greater than $v$ if and only if $\#B(x,2v)\cap \mathbf{X} = 1$. Consequently
\[\PPP{r(\cell) > v} = \PP^0(\# B(0,2v)\cap \mathbf{X}^0 = 1)\] where $\cell$ is the typical cell of the Voronoi tessellation induced by $\mathbf{X}$. In the above equality, $\PP^0$ is the Palm measure of $\mathbf{X}$ in the sense of (3.6) of \cite{SW} and $\mathbf{X}^0$ is $\PP^0$ distributed. The planar Gauss-Poisson process is one of the rare non-Poisson processes for which the right-hand side can be made fully explicit. This one is given for each $v\geq 0$ by (see p. 161 in \cite{SKM}):
 \begin{equation}
 \label{GaussPoisstypicaldistance}
 \PP^0(\# B(0,2v)\cap \mathbf{X}^0 = 1) = \frac{1}{p_1+2p_2}e^{-\gamma_a(4p_1\pi v^2+p_2(8\pi v^2 - a(2v)))}\cdot \left\{\begin{split} &  p_1+2p_2 \hspace{0.5cm} & 0\leq 2v<1\\ & p_1 & 2v\geq 1 \end{split}\right. .
 \end{equation} and \begin{equation}\label{defintersect}a(2v)=8v^2\arccos\frac{1}{4v} - \frac{1}{2}\sqrt{16v^2-1} \text{ for } 4v\geq 1\end{equation} and equals zero otherwise. The function $a(2v)$ is the area of intersection of two disks of radius $2v$ and centers separated by unit distance. A Taylor expansion of the right-hand side of \eqref{GaussPoisstypicaldistance} shows that 
 \begin{equation}\label{defPRGP}\PP^0(\# B(0,2v)\cap \mathbf{X}^0 = 1) = e^{-(P(v)+R(v))}\end{equation} where
 \begin{equation} \label{defPR}P(v) = 4\gamma_a\pi (p_1+p_2)v^2 - 4\gamma_a\cdot p_2\cdot v -\log\left(\frac{p_1}{p_1+p_2} \right) \text{ and } R(v) = \frac{5\gamma_a\cdot p_2}{48}\cdot \frac{1}{v}+o\left(\frac{1}{v} \right)\end{equation} as $v$ goes to infinity. In the previous line, $\phi(v) = o(\psi(v))$ means that $\phi(v)/\psi(v)\conv[v]{\infty}0$. 
 
For all $t\in \RR$, we define $v_{\rho}=v_{\rho}(t)$ so that $P(v_{\rho}) = \log\rho + t$ i.e. 
\begin{equation}\label{defvGP}v_{\rho} = v_{\rho}(t) = \frac{2\gamma_a\cdot p_2 + \left(4\gamma_a^2\cdot p_2^2 + 4\gamma_a\pi (p_1+p_2) \left(\log\left(\frac{p_1}{p_1+2p_2} \right)  + \log\rho + t \right)\right)^{1/2}}{4\gamma_a\pi (p_1+p_2)}.\end{equation} Using the fact that $\rho\PP^0(\# B(0,2v_{\rho})\cap \mathbf{X}^0 = 1) = e^{-t-R(v_{\rho})}$ where $R(\cdot)$ is defined in \eqref{defPR}, we deduce that
\begin{equation}
\label{typicalGP}
G_1(\rho) = |\rho\PP^0(\# B(0,2v)\cap \mathbf{X}^0 = 1) - e^{-t}| \leq e^{-t}R(v_{\rho}) = O\left((\log\rho)^{-1/2}\right).
\end{equation}

 Moreover, from Campbell theorem (see Theorem 3.3.3. in \cite{SW}), we have 
\begin{multline*}
G_2(\rho) := N_{\rho}\EEE{\sum_{(x,y)_{\neq} \in (\mathbf{X}\cap \CC)^2}\mathbb{1}_{\# B(x,2v_{\rho})\cap \mathbf{X} =1}\mathbb{1}_{\# B(y,2v_{\rho})\cap \mathbf{X} =1}}\\
 = N_{\rho}\int_{\mathfrak{C}_{\rho}}\int_{\mathcal{F}_{lf}}\sum_{y\in \eta\cap \CC}\mathbb{1}_{\#(\eta +x)\cap B(x,2v_{\rho}) = 1 }\mathbb{1}_{\#(\eta +x)\cap B(y,2v_{\rho}) = 1 }d\PP^0(\eta)dx.
\end{multline*} 
Here $\mathcal{F}_{lf}$ denotes the space of locally finite subsets of $\RR^2$. Because the integrand of the right-hand side is translation invariant (in distribution) and because $N_{\rho}\lambda_2(\mathfrak{C}_{\rho}) = c\cdot \rho$, we obtain
\begin{equation*}
G_2(\rho) = c\cdot\rho\int_{\mathcal{F}_{lf}}\sum_{y\in \eta\cap \CC}\mathbb{1}_{\#\eta\cap B(0,2v_{\rho}) = 1 }\mathbb{1}_{\#\eta\cap B(y,2v_{\rho}) = 1 }d\PP^0(\eta).
\end{equation*}
 According to Formula (5.3.2) in \cite{SKM}, we have $\PP^0 = \PP_{\mathbf{X}}\ast \mathbb{c}^0$ where $\PP_{\mathbf{X}}$ is the distribution of $\mathbf{X}$ and $\mathbb{c}^0$ is the Palm measure of the  cluster distribution $\Xi_0$ that is concentrated on the space $\mathcal{F}_{lf,2}$ of subsets of 0, 1 or 2 points in $\RR^2$. Hence 
\[G_2(\rho) = c\cdot \rho \int_{\mathcal{F}_{lf}}\int_{\mathcal{F}_{lf,2}}\sum_{y\in (\phi\cup \xi)\cap \CC}   \mathbb{1}_{\#(\phi\cup \xi)\cap (B(0,2v_{\rho})\cup B(y,2v_{\rho})) = 2}\mathbb{1}_{|y|>2v_{\rho}}d\PP_{\mathbf{X}}(\phi)d\mathbb{c}^0(\xi).\] When $|y|>2v_{\rho}$, we have $y\not\in \xi$ for $\rho$ large enough since $\mathbb{c}_0$ a.s. $\xi$ is bounded. Moreover, $\PP_{\mathbf{X}}$ a.s. $\phi\cap \xi\cap (B(0,2v_{\rho})\cup B(y,2v_{\rho}))$ is empty. Consequently, calculating the integral with respect to $\mathbb{c}_0$, we get
\begin{equation*}G_2(\rho) = c\cdot \rho \int_{\mathcal{F}_{lf}}\sum_{y\in \phi\cap \CC}\mathbb{1}_{\#\phi\cap (B(0,2v_{\rho})\cup B(y,2v_{\rho})) = 1}\mathbb{1}_{|y|>2v_{\rho}}d\PP_{\mathbf{X}}(\phi).\end{equation*} 
Proceeding as previously, we deduce from Campbell theorem and from the relation $\PP^0 = \PP_{\mathbf{X}}\ast \mathbb{c}^0$, that
\begin{equation*}
G_2(\rho) = c\cdot \rho\int_{\CC}\int_{\mathcal{F}_{lf}}\int_{\mathcal{F}_{lf,2}}\mathbb{1}_{\#((\xi\cup\phi)+y)\cap (B(0,2v_{\rho})\cup B(y,2v_{\rho})) = 1}\mathbb{1}_{|y|>2v_{\rho}}dyd\PP_{\mathbf{X}}(\phi)d\mathbb{c}^0(\xi).
\end{equation*}Since $\PP_{\mathbf{X}}$ a.s. $\phi\cap \Xi_0\cap (B(0,2v_{\rho})\cup B(y,2v_{\rho}))$  is empty, we deduce after integration over $\mathcal{F}_{lf}\times \mathcal{F}_{lf,2}$ with respect to $\PP^0\otimes\mathbb{c}^0$ that \begin{equation}\label{GP1}G_2(\rho) \leq c\cdot \rho\int_{\CC}\PPP{\mathbf{X}\cap (B(0,2v_{\rho})\cup B(y,2v_{\rho})) = \varnothing} \mathbb{1}_{|y|>2v_{\rho}}dy.\end{equation} Let $|y|>2v_{\rho}$ be fixed. To get a suitable upper bound of the integrand, we use the fact that $\mathbf{X}\cap (B(0,2v_{\rho})\cup B(y,2v_{\rho}))=\varnothing\Longleftrightarrow (x+\Xi_0(x))\cap (B(0,2v_{\rho})\cup B(y,2v_{\rho})) = \varnothing$ for all $x\in \mathbf{X}_a$. From Theorem 3.2.4. of \cite{SW}, Fubini's theorem and the fact that $\Xi_0$ is symmetric, we get
 \begin{equation}\label{GP2}\begin{split} \PPP{\mathbf{X}\cap (B(0,2v_{\rho})\cup B(y,2v_{\rho})) = \varnothing} & = e^{-\gamma_a\int_{\RR^2}\PPP{(x+\Xi_0(x))\cap (B(0,2v_{\rho})\cup B(y,2v_{\rho}))\neq\varnothing}dx}\\
 & = e^{-\gamma_a\EEE{\lambda_2(\Xi_0\oplus (B(0,2v_{\rho})\cup B(y,2v_{\rho})))}}.\end{split}\end{equation} 
We give below a suitable lower bound of the term appearing in the exponential. Since $|y|>2v_{\rho}$, we have
\[\EEE{\lambda_2(\Xi_0\oplus (B(0,2v_{\rho})\cup B(y,2v_{\rho})))|\#\Xi_0 = 1} = \lambda_2\left(B(0,2v_{\rho})\cup B(y,2v_{\rho}) \right)\geq \frac{3}{2}\cdot 4\pi v_{\rho}^2 \]  and
\[\EEE{\lambda_2(\Xi_0\oplus (B(0,2v_{\rho})\cup B(y,2v_{\rho})))|\#\Xi_0 = 2} \geq \EEE{\lambda_2(\Xi_0\oplus B(0,2v_{\rho}))}\geq  8\pi v_{\rho}^2 - a(2v_{\rho})\] where  $a(\cdot)$ is defined in \eqref{defintersect}. Since $\Xi_0$ is reduced to 0, 1 or 2 points with probability $p_0$, $p_1$ and $p_2$, we deduce from \eqref{GP2} that
\begin{multline}\label{GP3}\PPP{\mathbf{X}\cap (B(0,2v_{\rho})\cup B(y,2v_{\rho})) = \varnothing}\leq e^{-\gamma_a\left(\frac{3}{2}p_1\cdot 4\pi v_{\rho}^2 + p_2(8\pi v_{\rho}^2 - a(2v_{\rho})) \right)  }\\  = \frac{p_1+2p_2}{p_1} \PP^0(\# B(0,2v_{\rho})\cap \mathbf{X}^0 = 1)\cdot e^{-2\gamma_ap_1\pi v_{\rho}^2}   \end{multline} for $\rho$ large enough according to \eqref{GaussPoisstypicaldistance}. Integrating over $\mathfrak{C}_{\rho}$, we deduce from \eqref{typicalGP}, \eqref{GP1}, \eqref{GP3}  and from the inequality $\lambda_2(\mathfrak{C}_{\rho})\leq c\cdot \log\rho$, that
\begin{equation}\label{order2GP}G_2(\rho) \leq c\cdot \log\rho\cdot e^{-2\gamma_ap_1\pi v_{\rho}^2} = O\left(\log\rho\cdot \rho^{-\alpha} \right)\end{equation} where \begin{equation}\label{defalphaconv}\alpha = \frac{p_1}{2(p_1+p_2)}.\end{equation} Since $p_1\neq 0$, we have $\alpha>0$ so that $G_2(\rho)$ converges to 0. Proposition \ref{GaussPoissonprop} is now a direct consequence of \eqref{typicalGP}, \eqref{order2GP} and Theorem \ref{realtyp}.

\end{prooft}

According to Proposition \ref{GaussPoissonprop} and \eqref{defvGP}, the order of $r_{\max,GPVT}(\rho)$ is 
\[\left(4\gamma_a\pi (p_1+p_2)\right)^{-1/2}\cdot (\log\rho)^{1/2} = \left(\frac{p_1+2p_2}{4\pi(p_1+p_2)} \right)^{1/2}\cdot (\log\rho)^{1/2}\] since we have assumed that $(p_1+2p_2)\gamma_a=1$. Let us remark that the larger $p_2$ is, the larger the order is. This can be explained by the following fact: the nucleus $x\in \mathbf{X}$ of the Voronoi cell which maximizes the inradius belongs to a cluster of size 1 i.e. $x\in \Xi(x_a)$ where $\#\Xi(x_a)=1$ for some $x_a\in \mathbf{X}_a$. Hence if $p_2$ is large, the mean number of clusters of size 1 is small so that the inradii associated to the clusters of size 1 are large.

When $p_1=0$, we obtain a degenerate case since $r_{\max, GPVT}(\rho)=\frac{1}{2}$ is constant. When $p_0=p_2=0$ and $p_1=1$, the random variable $r_{\max,GPVT}(\rho)$ is the maximum of inradii of a Poisson-Voronoi tessellation $r_{\max,PVT}(\rho)$. In that case, the order is
\[v_{\rho} = v_{\rho}(t) = (4\pi)^{-1/2}\cdot\left(\log\rho + t\right)^{1/2}.\] The  order  of $r_{\max,PVT}(\rho)$ has already been investigated in \cite{CC}. Nevertheless, Proposition \ref{GaussPoissonprop} is more precise since it provides the rate of convergence. Actually, this rate could be improved. Indeed, since $p_0=p_2=0$ and $p_1=1$ we have $R(v_{\rho}) = 0$ according to \eqref{GaussPoisstypicaldistance} and \eqref{defPRGP} and consequently we get $G_1(\rho)=0$ according to the inequality in \eqref{typicalGP}. Moreover, the term $\alpha$ as defined in \eqref{defalphaconv} equals $1/2$. Hence, according to \eqref{order2GP}, we obtain the more precise result:
\[\PPP{r_{\max, PVT}(\rho) \leq (4\pi)^{-1/2}\cdot\left(\log\rho + t\right)^{1/2}} = O\left(\log\rho\cdot \rho^{-1/2} \right).\]

Finally, let us mention that a Gauss-Poisson process belongs to the class of the so called Neyman-Scott processes. We do not investigate general Neyman-Scott processes since the left-hand side of \eqref{GaussPoisstypicaldistance} cannot be made explicit excepted for some particular cases as Gauss-Poisson processes.

\section{Proof of Proposition \ref{realtyp2} and some extremal indices}
\label{sectionrealtyp2}
In this section, we prove Proposition \ref{realtyp2} and we give two examples where the extremal index differs from 1.   

\begin{prooft}{Proposition \ref{realtyp2}}
The proof is an adaptive version to our setting of two results due to Leadbetter (see Theorem 2.2 and Lemma 2.1. in \cite{L2}). The difference is that we investigate a maximum on a random graph instead of a sequence of real numbers. 

First, we investigate the limit superior. For each $\tau\geq 0$, we denote by 
\begin{equation}
\label{defpsi}
\psi(\tau) = \limsup_{\rho\rightarrow\infty}\PPP{M_{f,\mathbf{W}_{\rho}}\leq v_{\rho}(\tau)}.
\end{equation}
Let $\tau\geq 0$ and $k\in \NN^*$ be fixed. The key idea is to show that $\psi(\tau/k^d)=\psi^{1/k^d}(\tau)$. To do it, we subdivide the proof into two steps. The first is intrinsic to the sequence $v_{\rho}(\tau)$ while the second step needs the mixing property of the tessellation i.e. \textsc{Condition 1}. 

\textbf{Step 1.} We show that \begin{equation}
\label{step1prop2}
\limsup_{\rho\rightarrow\infty}\PPP{M_{f,\mathbf{W}_{\rho/k^d}}\leq v_{\rho}(\tau)} = \psi(\tau/k^d).
\end{equation}
Indeed, if $v_{\rho}(\tau)\geq v_{\rho/k^d}(\tau/k^d)$, it follows that \begin{multline*}
\left|\PPP{M_{f,\mathbf{W}_{\rho/k^d}}\leq v_{\rho}(\tau)} - \PPP{M_{f,\mathbf{W}_{\rho/k^d}}\leq v_{\rho/k^d}(\tau/k^d)} \right| \leq \PPP{\bigcup_{\underset{z(C)\in \mathbf{W}_{\rho/k^d}}{C\in\mm}}\{v_{\rho/k^d}(\tau/k^d) \leq f(C)\leq v_{\rho}(\tau)\}}\\
\leq \EEE{\sum_{\underset{z(C)\in \mathbf{W}_{\rho/k^d}}{C\in\mm}}\mathbb{1}_{v_{\rho/k^d}(\tau/k^d) \leq f(C)\leq v_{\rho}(\tau)}}.
\end{multline*} This together with the corresponding inequality when $v_{\rho}(\tau)\leq v_{\rho/k^d}(\tau/k^d)$ shows that
\begin{multline}
\label{step11prop2}
\left|\PPP{M_{f,\mathbf{W}_{\rho/k^d}}\leq v_{\rho}(\tau)} - \PPP{M_{f,\mathbf{W}_{\rho/k^d}}\leq v_{\rho/k^d}(\tau/k^d)} \right| \leq \frac{\rho}{k^d}\left|\PPP{f(\cell)>v_{\rho/k^d}(\tau/k^d)} - \PPP{f(\cell)>v_{\rho}(\tau)} \right|\\
 = \frac{\rho}{k^d}\left|\frac{\tau/k^d}{\rho/k^d} - \frac{\tau}{\rho}+o\left(\frac{1}{\rho} \right) \right|\conv[\rho]{\infty}0
\end{multline} according to \eqref{campbell} and the fact that $\PPP{f(\cell)>v_{\rho}(\tau)}$ converges to $\tau$ for each $\tau\geq 0$ . Moreover, from \eqref{defpsi} we have \[\limsup_{\rho\rightarrow\infty}\PPP{M_{f,\mathbf{W}_{\rho/k^d}}\leq v_{\rho/k^d}(\tau/k^d)} = \psi(\tau/k^d).\] The limit \eqref{step1prop2} results of the previous equality and \eqref{step11prop2}. 

\textbf{Step 2.} Secondly, we show that \begin{equation}\label{step2prop2}
\limsup_{\rho\rightarrow\infty}\PPP{M_{f,\mathbf{W}_{\rho/k^d}}\leq v_{\rho}(\tau)} = \psi(\tau)^{1/k^d}.\end{equation}
Indeed, we partition $W=[0,1]^d$ into a set of $k^d$ sub-cubes of equal volume $1/k^{d}$ say $B^{(1)}, \ldots, B^{(k^d)}$. According to \eqref{Rkindependence} applied to $L=k^d$, we have
\[\PPP{M_{f,\mathbf{W}_{\rho}}\leq v_{\rho}(\tau)} - \prod_{l=1}^{k^d}\PPP{M_{f,\mathbf{B}_{\rho}^{(l)}}\leq v_{\rho}(\tau)}\conv[\rho]{\infty}0\] where $\mathbf{B}_{\rho}^{(l)} = \rho^{1/d}B^{(l)}$ for all $1\leq l\leq k^d$. Since $\mathbf{B}_{\rho}^{(l)}$ is a cube of volume $\rho/k^d$ and since $\mm$ is stationary, the previous convergence can be re-written as \begin{equation}\label{independentproba}\PPP{M_{f,\mathbf{W}_{\rho}}\leq v_{\rho}(\tau)} - \PPP{M_{f,\mathbf{W}_{\rho/k^d}}\leq v_{\rho}(\tau)}^{k^d}\conv[\rho]{\infty}0.\end{equation} We deduce \eqref{step2prop2} thanks to \eqref{defpsi}.

\textbf{Conclusion}. We deduce from \eqref{step1prop2} and \eqref{step2prop2}, that \begin{equation}\label{functionalequation}\psi(\tau/k^d) = \psi(\tau)^{1/k^d} \text{ where } \tau\geq 0 \text{ and } k\in\NN^*  \end{equation} Moreover, in the same spirit as in the proof of \eqref{majcubei}, we can show that
\[\PPP{ M_{f, \mathbf{W}_{\rho/k^d}}\leq v_{\rho}(\tau) }\geq 1-\frac{\rho}{k^d}\PPP{f(\cell)>v_{\rho}(\tau)}\conv[\rho]{\infty}1-\tau/k^d.\] Hence, taking the $k^{\text{th}}$ powers and using \eqref{independentproba}, we deduce that $\PPP{M_{f, \mathbf{W}_{\rho}}\leq v_{\rho}(\tau)} \geq \left(1-\frac{\tau}{k^d}\right)^{k^d}$ and so, letting $k\rightarrow\infty$, that
\begin{equation}\label{positivefunction}\liminf_{\rho\rightarrow\infty}\PPP{M_{f, \mathbf{W}_{\rho}}\leq v_{\rho}(\tau)}\geq e^{-\tau}. \end{equation} This shows that $\psi(\tau)>0$. Since $\psi(\cdot)$ is also non-increasing and since the only solution of the functional equation \eqref{functionalequation} which is strictly positive and non-increasing is an exponential function, we have $\psi(\tau)=e^{-\theta\tau}$ for some $\theta\geq 0$. Hence
\[\limsup_{\rho\rightarrow\infty }\PPP{M_{f, \mathbf{W}_{\rho}}\leq v_{\rho}(\tau)} = e^{-\theta\tau}.\]

With a similar method, we obtain that $\liminf_{\rho\rightarrow\infty }\PPP{M_{f, \mathbf{W}_{\rho}}\leq v_{\rho}(\tau)} = e^{-\theta'\tau}$ for some $\theta'\leq 1$ (according to \eqref{positivefunction}) and such that $\theta\leq \theta'$.
 \end{prooft}
 
As an illustration, we give below two examples where the extremal index differs from 1. The first one is the minimum of inradii of a Poisson-Voronoi tessellation.

\begin{Ex}
Let $\mm_{PVT}$ be a Poisson-Voronoi tessellation of intensity 1 and $\mathbf{X}$ the underlying Poisson point process. For each cell $C=C_{\mathbf{X}}(x)$, we consider the inradius $r(C_{\mathbf{X}}(x))$ in the sense of section \ref{sectionGP} and we denote by $r(\cell)$ the inradius of the typical cell $\cell \overset{\mathcal{D}}{=} C_{\mathbf{X}\cup\{0\}}(0)$. The distribution function of $r(\cell)^d$  is exponentially distributed with rate $2^d\kappa_d$. Indeed, $r(\cell)$ is lower than $v$, $v\geq 0$, if and only if $\mathbf{X}\cap B(0,2v)$ is not empty. Hence
 \begin{equation*}
\label{mininstyp}
\rho\cdot\PPP{r(\cell)^d \leq \frac{1}{2^d\kappa_d\rho}t}\conv[\rho]{\infty}t.
\end{equation*} Moreover, according to the convergence (2b) in \cite{CC}, we know that
\begin{equation*}
\PPP{\min_{x\in \mathbf{X}\cap \mathbf{W}_{\rho}}r(C_{\mathbf{X}}(x))^d \geq \frac{1}{2^d\kappa_d\rho}t}\conv[\rho]{\infty}e^{-t/2}.
\end{equation*}
Let us notice that the convergence was written in \cite{CC} for a fixed window and for a Poisson point process such that the intensity goes to infinity. By scaling property of the Poisson point process, the result of \cite{CC} can be re-written as above for a fixed intensity and for a window $W_{\rho}$ where $\rho\rightarrow\infty$.
 
Therefore, the extremal index of the minimum of inradii is \[\theta=\frac{1}{2}.\]
It can be aslo explained by a trivial heuristic argument. Indeed, if a cell minimizes the inradius, one of its neighbors has to do the same. Hence, the mean cluster size of exceedances is 2. This justifies the fact that $\theta=1/2$. 
\end{Ex}

In our second example, we give the extremal index of the maximum of circumradii of a Poisson-Delaunay tessellation. 

\begin{Ex}
Let $\mm_{PDT}$ be a Poisson-Delaunay tessellation of intensity 1 and let $\mathbf{X}$ be the underlying Poisson point process (of intensity $\gamma_{\mathbf{X}} = \beta_d^{-1}$ where $\beta_d^{-1}$ is given in \eqref{defad}). Denoting by $\cell$ the typical cell of $\mm_{PDT}$, we deduce from a Taylor expansion of \eqref{circonsdel} that 
\begin{equation*}
\rho\cdot\PPP{R(\cell)^d \geq \frac{\log\left([(d-1)!]^{-1}\rho\log(\beta_d\rho)^{d-1} \right) + t}{\delta_d}}\conv[\rho]{\infty}e^{-t}
\end{equation*} for all $t\in\RR$. Moreover, considering the dual Voronoi tessellation of $\mm_{PDT}$, we have
\begin{equation}
\label{maxcirconsdelvor}
\max_{x\in \mathbf{X}\cap \mathbf{W}_{\rho}}R(C_{\mathbf{X}}(x)) = \max_{\underset{V(C)\cap \mathbf{W}_{\rho}\neq\varnothing}{C\in \mm_{PDT}}}R(C) 
\end{equation} where $V(C)$ is the set of vertices of the Delaunay cell $C\in\mm_{PDT}$. The asymptotic behaviour of the maximum of circumradii of a Poisson-Voronoi tessellation is already known (see (2c) in \cite{CC}). This is given by 
\begin{equation}\label{maxcirconsVor}\PPP{\max_{x\in \mathbf{X}\cap \mathbf{W}_{\rho}}R(C_{\mathbf{X}}(x))^d \leq \frac{\log\left(\alpha_{d,6}\beta_d\rho\log(\beta_d\rho)^{d-1}\right) + t}{\delta_d}}\conv[\rho]{\infty}e^{-e^{-t}}\end{equation} where $\alpha_{d,6}:=\frac{1}{d!}\left(\frac{\pi^{1/2}\Gamma\left(\frac{d}{2}+1\right)}{\Gamma\left(\frac{d+1}{2}\right)}\right)^{d-1}$. With a similar method as in Lemma 4.1. in \cite{HSS}, we can show that the boundary cells of the Poisson-Delaunay tessellation (i.e. the cells which intersect the boundary of $\mathbf{W}_{\rho}$) do not affect the behaviour of the maximum. Hence, the rate of $\max_{\underset{V(C)\cap \mathbf{W}_{\rho}\neq\varnothing}{C\in \mm_{PDT}}}R(C) $ is the same as $\max_{\underset{z(C)\in \mathbf{W}_{\rho}}{C\in\mm_{PDT}}}R(C)$. We then deduce from  \eqref{maxcirconsdelvor} and \eqref{maxcirconsVor} that \begin{equation}\label{maxcircumdel}\PPP{\max_{\underset{z(C)\in \mathbf{W}_{\rho}}{C\in\mm_{PDT}}}R(C)^d \leq \frac{\log\left([(d-1)!]^{-1}\rho\log(\beta_d\rho)^{d-1} \right) + t}{\delta_d} }\conv[\rho]{\infty} e^{-e^{-t}\times\theta}\end{equation} where \[\theta = \alpha_{d,6}\beta_d(d-1)! = \frac{(d^3+d^2)\Gamma\left(\frac{d^2}{2} \right) \Gamma\left(\frac{d+1}{2}\right)}{2^{d+1}d\Gamma\left(\frac{d^2+1}{2} \right)\Gamma\left(\frac{d+2}{2} \right)}.\] In particular, when $d=1,2,3$, the extremal indices are $\theta=1$, $\theta=1/2$ and $\theta=35/128$ respectively. The fact that $\theta=1$ when $d=1$ follows from Theorem \ref{realtyp} which is available since the associated function $G_2(\cdot)$ converges to 0. This is not the case in higher dimension.
\end{Ex}

We hope to be able to develop a systematic method to estimate the extremal index in a future work.

\paragraph{Acknowledgements.} I would like to thank my advisor P. Calka for helpful discussions and suggestions. This work was partially supported by the French ANR grant PRESAGE (ANR-11-BS02-003) and the French research group GeoSto (CNRS-GDR3477).


\bibliographystyle{abbrv}
\bibliography{BiblioVE.bbl}

\end{document}